\documentclass{article} \usepackage{amssymb,amsmath,eucal}

\begin{document}

\newcounter{sec} \renewcommand{\theequation}{\arabic{sec}.\arabic{equation}}


\newcounter{punct}
\def\punct{\addtocounter{punct}{1}{\arabic{sec}.\arabic{punct}.  }}

\newtheorem{theorem}{Theorem}[sec]
\newtheorem{prop}[theorem]{Proposition}
\newtheorem{lemma}[theorem]{Lemma}
\newtheorem{cor}[theorem]{Corollary}

\def\R {{\mathbb R }}
 \def\C {{\mathbb C }}
  \def\Z{{\mathbb Z}}
  \def\H{{\mathbb H}}
\def\K{{\mathbb K}}

\def\T{\mathbb T}

\newcommand{\arcsh}{\mathop{\rm arcsh}\nolimits}
 \newcommand{\sh}{\mathop{\rm sh}\nolimits}

\newcommand{\SU}{\mathop{\mathrm{SU}}\nolimits}

\def\Pol{\mathrm{Sym}}
\def\LPol{\mathrm{LSym}}

\def\SL{\mathrm {SL}}
\def\GL{\mathrm  {GL}}
\def\U{\mathrm  U}
\def\OO{\mathrm  O}
\def\Sp{\mathrm  {Sp}}
\def\SO{\mathrm  {SO}}
\def\SOS{\mathrm {SO}^*}

\def\Mat{{\mathrm {Mat}}}
\def\Symm{{\mathrm {Symm}}}
\def\ASymm{{\mathrm {ASymm}}}
\def\Herm{{\mathrm {Herm}}}
\def\AHerm{{\mathrm {AHerm}}}
\def\Sti{\mathrm Sti}

\def\cY{{\mathcal Y}}
\def\cX{{\mathcal X}}
\def\cM{{\mathcal M}}
\def\cD{{\mathcal D}}
\def\cH{{\mathcal H}}
\def\cB{{\mathcal B}}
\def\cK{\mathcal K}
\def\cL{\mathcal L}
\def\cA{\mathcal A}
\def\cB{\mathcal B}
\def\cC{\mathcal C}

\def\frD{\mathfrak D}

 \def\ov{\overline}
\def\wt{\widetilde}

\def\phi{\varphi}
\def\epsilon{\varepsilon}
\def\kappa{\varkappa}

\def\le{\leqslant}
\def\ge{\geqslant}

\def\R{{\mathbb R}}
\def\C{{\mathbb C}}
\def\H{{\mathbb H}}
\def\K{{\mathbb K}}

\def\F{\,\,{}_2F_1}
\def\FF{\,\,{}_3F_2}

\renewcommand{\Re}{\mathop{\rm Re}\nolimits}
\renewcommand{\Im}{\mathop{\rm Im}\nolimits}

\begin{center} \Large

\bf Notes on  Stein-Sahi representations
and some problems of non $L^2$ harmonic analysis

\medskip

\large
\sc Neretin Yu.A.

\end{center}

{\small
We discuss one natural class of kernels
on pseudo-Riemannian symmetric spaces.}

\medskip
Recently, Oshima \cite{Osh} published his formula
for $c$-function for
$L^2$ on pseudo-Riemannian symmetric spaces
(see also works of Delorm \cite{Del}
and van den Ban--Schlichtkrull \cite{Ban1}, \cite{Ban2}).
After this, there arises a natural question
about other solvable problems of
non-commutative  harmonic analysis.

In the Appendix to the paper
\cite{Ner-determinant},
the author proposed a series of 
non $L^2$-inner products in spaces of functions
on pseudo-Riemannian symmetric spaces
and conjectured that this object is reasonable and admits
an explicit harmonic analysis.

In this work, we discuus the problem in more details,
 in particular, we obtain the Plancherel formula
 for these kernel for Riemannian
 symmetric spaces $\U(n)$, $\U(n)/\OO(n)$,
 $\U(2n)/\Sp(n)$.
 We also give a new proof of Sahi's results
 \cite{Sah2}.

\bigskip

{\bf\Large 0. Introduction}

\bigskip

{\bf \punct Inner products defined by kernels.} Starting the
famous works of Bargmann \cite{Bar} and Gelfand--Naimark
\cite{GN1}, various inner products having the form
\begin{equation}
\langle f_1, f_2 \rangle=\iint_{G/H\times G/H} K(x,y)
f_1(x)\ov{f_2(y)}dx\,dy \label{yadro}
\end{equation}
 are quite usual in the  theory of unitary
representations. Here $G/H$ is a homogeneous space and $K(x,y)$ is
a distribution (a {\it 'kernel'}) on $G/H\times G/H$. The group
$G$ acts in a space of functions on $G/H$ by transformations
having the form
\begin{equation}
\rho(g)f(x) = f(gx)\gamma(g,x) \label{action}
,\end{equation}
where $\gamma(g,x)$ is some function ({\it
'multiplier'}) on $G\times G/H$. We intend to obtain a unitary
irreducible
representation; under this requirement,  the kernel $K(x,y)$
 is uniquely determined
 by the
explicit expression for the multiplier $\gamma(g,x)$. An actual
evaluation of the kernel  is not difficult.

But the scalar square $\langle f,f\rangle$ of a function $f$,
i.e., the integral
$$
\langle f, f \rangle=\iint_{G/H\times G/H} K(x,y)
f(x)\ov{f(y)}dx\,dy
$$
  has no visible
reasons to be positive. Usually,  positive definiteness of a given
inner product of the form (\ref{yadro}) is a nontrivial problem.

\smallskip

{\sc Example.} We consider group $\SU(1,1)$
consisting of $2\\times 2$ matrices
$\begin{pmatrix}a&b\\\ov b & \ov a\end{pmatrix}$,
where $|a^2|-|b|^2=1$.
It acts in the space of functions
on the circle $|z|=1$ (or $z=e^{iphi}$) by the formula
\begin{equation}
\rho_s  
\begin{pmatrix}a&b\\\ov b & \ov a\end{pmatrix}
f(z)=f\Bigl(\frac{az+b}{\ov b z+\ov a}\Bigr)
|\ov b z+\ov a|^{1+s}
\label{sl2}
\end{equation}
These operators preserve the inner product given by
\begin{equation}
\langle f,g\rangle
=\iint
\bigl|z_1-z_2|^{-1-s}f_1(z_1)\ov {f_2(z_2)}
\,d\phi_1\, d\phi_2
\label{sobolev-circle}
\end{equation}
This inner product is positive definite
iff $-1<s<1$.

\smallskip

The inner product (\ref{sobolev-circle})
is a Sobolev inner product. Also, the Sobolev inner products
in spaces of functions on spheres
appears in the representation theory of $\SO(1,n)$.
Some anisotropic Sobolev spaces arise
in the representation
theory of rank 1 groups. But usually inner products
(\ref{yadro}) define functional spaces that
are uknown in  the analysis. The spaces
discussed in this paper can be concidered
 as some kind
of "{\it Sobolev spaces of matrix variables}.

\smallskip

{\bf \punct Stein--Sahi complementary series.} In 1967, E. Stein
\cite{Ste} constructed an extremely degenerated complementary
series of unitary   representations of $\GL(2n,\C)$. Recall his
construction. Consider the space $\Mat_n(\C)$ consisting  of
complex $n\times n$ matrices. Consider the group $\GL(2n,\C)$
consisting of $(n+n)\times (n+n)$ invertible complex matrices
$g=\begin{pmatrix}a&b\\c&d\end{pmatrix}$. This group acts on
$\Mat_n(\C)$ by linear-fractional transformations
\begin{equation}
z\mapsto (a+zc)^{-1}(b+zd) \label{lin-fra}
.\end{equation}
Fix $\sigma\in\R$
Consider the action of $\GL(2n,\C)$ in the space of
functions on $\Mat_n(\C)$ by the operators
$$
\rho_\sigma(g)\,f(z)= f\bigl( (a+zc)^{-1}(b+zd) \bigr)
|\det(a+zc)|^{-2n-\sigma}
.$$
Define the following Hermitian form in the space
$\cD(\Mat_n(\C))$ of smooth functions  on $\Mat_n(\C)$ with
compact support
\begin{equation}
\langle f,g\rangle_\sigma=
\iint_{\Mat_n(\C)\times\Mat_n(\C)}
|\det(x-y)|^{-2n+\sigma} f(x)\ov {g(y)}\,dx\,dy
.
\label{STEIN}
\end{equation}
For $\sigma>2n$ this integral converges, further we consider its
meromorphic  continuation in $\sigma$ to the whole plane
$\sigma\in\C$.

Stein proved that for $-1<\sigma<1$,

\smallskip

 1. the Hermitian form $\langle\cdot,\cdot\rangle_\sigma$
is positive definite; denote by
$\cH_\sigma$ the completion of $\cD(\Mat_n(\C))$ with respect to
$\langle\cdot,\cdot\rangle_\sigma$

 2. the operators $\rho_\sigma(g)$ are
 unitary in $\cH_\sigma$.

\smallskip

The Stein-type constructions exist for all the series
of classical groups, precisely for the groups
\begin{align}
\GL(2n,\C),\, \GL(2n,\R),\, \GL(2n,\R),\, \label{line1}
\\
\OO(2n,2n),\, \U(n,n),\, \Sp(n,n),\label{line2}
\\
\SO(4n,\C),\, \Sp(2n,\C),\, \SOS(4n,\C),\,\Sp(2n,\R)\label{line3}
.\end{align}

The cases (\ref{line1}) were considered in the
classification work of
Vogan (\cite{Vog}, Section 2), see also  \cite{Sah}.
Sahi \cite{Sah2}--\cite{Sah3}
 constructed  analogs of the Stein
representations in the remaining cases.

\smallskip

\smallskip

{\bf \punct Berezin kernels and problem of existence of
non-$L^2$-theories.} The explicit Plancherel formula for $L^2$ on
Riemannian noncompact symmetric spaces was obtained by Gindikin
and Karpelevich  \cite{GK2} in 1964. In 1978, Berezin \cite{Ber2}
       \footnote{   The first work that can be attributed to this
       subject, is Vershik,
   Gelfand, Graev, \cite{VGG}, \cite{GGV}. These authors
 apply the representations $\rho_\sigma$ of the groups $\U(1,q)$ to
       construct representations of current groups;
for a collection of other early references, see \cite {Ner-berezin}.}
  observed that the space $L^2$ on a
  classical  Hermitian symmetric
space admits a natural  deformation.

For definiteness,
consider the Hermitian symmetric space
$$G/K=\U(p,q)/\U(p)\times\U(q), \qquad p\le q$$
Following E.Cartan, we realize this space as
 the space
$\cB_{p,q}$ of $p\times q$ complex matrices with norm $<1$, the
group $\U(p,q)$ acts on this space by the linear-fractional
transformations $ z\mapsto (a+zc)^{-1} (b+zd) $.
 Consider (see
\cite{Ner-determinant} for details)
 the space $\cD(\cB_{p,q})$
of smooth function on $\cB_{p,q}$ with compact support. Denote the
representation of the group $\U(p,q)$ in this space by the shifts
\begin{equation}
\rho_\sigma
\begin{pmatrix} a&b\\c&d\end{pmatrix}
F(z)= F\bigl((a+zc)^{-1} (b+zd)\bigr) \label{action-on-ball}
.\end{equation}
 For $\sigma\in\R$,
 define  the
inner product in $\cD(\cB_{p,q})$  by
\begin{equation}
\langle F_1, F_2\rangle_\sigma =
 \iint\limits_{\cB_{p,q}\times\cB_{p,q}}
\frac{\det(1-zz^*)^{\sigma}\det(1-uu^*)^{\sigma}}
     {|\det(1-zu^*)|^{2\sigma}}
     F_1(z)\,\ov {F_2(u)}\,d\sigma(z)\,d\sigma(u)
   \label{berezin-form}
     ,\end{equation}
where
$$
d\sigma(z)=\det(1-zz^*)^{-p-q} \prod_{1\le k\le p, 1\le l\le
q}\bigl( \frac1{2i}dz_{k,l}\,d\ov z_{kl}\bigr)
$$
is the $\U(p,q)$-invariant measure on $\cB_{p,q}$. A simple
calculation shows that the form (\ref{berezin-form}) is
$\U(p,q)$-invariant.%
\footnote{Compare the Stein inner product
(\ref{STEIN}) and the Berezin inner product (\ref{berezin-form}).
In the first case, the kernel has a singularity on the
surface $\det(x-y)=0$. This surface itself is singular,
and its  most singular strate is the diagonal $x=y$.
On the contary, the Berezin kernel has no singularities
in $\cB_{p,q}\times \cB_{p,q}$ at all.
We can substitute to
(\ref{berezin-form}) arbitrary compactly supported distributions
unstead of functions $F_1$, $F_2$.}

It turns out to be that {\it for
$$
\sigma=0,1,\dots, p-1,\quad\text{ or $\sigma>p-1$},
$$
our inner product is non-negative definite.}

The latter  statement is a reformulation of the well-known theorem
(Berezin--Gindikin--Rossi--Vergne--Wallach) on unitarizability of
scalar highest weight representations (\cite{Ber}, \cite{VR},
 see
a relatively recent exposition in \cite{FK}).

Consider the completion $\cH_\sigma$ of the space
$\cD(\cB_{p,q})$ with respect to the inner product
(\ref{berezin-form}), and the
unitary representation (\ref{action-on-ball}) of $\U(p,q)$ in this space.%
\footnote{This representation  is equivalent to the tensor
product of a scalar highest weight representation of $\U(p,q)$ and
the complex conjugate representation;  it is natural to
consider a similar and  more general problem about
the tensor product
of a highest weight and lowest weight representation, 
in this case, the Plancherel formula was obtained
by  Zhang,
\cite{Zha}.}


At first sight, the representation $\rho_\sigma$ given by the
formula (\ref{action-on-ball}) does not include $\sigma$.
Really, for sufficiently large $\sigma$, all  the representations
$\rho_\sigma$ are equivalent to the standard representation
of $\U(p,q)$ in $L^2(\U(p,q)/\U(p)\times\U(q)$.
Moreover, the natural limit of
the spaces as $\sigma\to\infty$ is
$L^2\bigl[\U(p,q)/\U(p)\times\U(q)\bigr]$

But the Plancherel formula for $\rho_\sigma$ depend on
$\sigma$.
For sufficiently small $\sigma>0$,
the spectrum of
 $\rho_\sigma$ undergo several
bifurcations.%
\footnote{In the work of Olshanski and the author \cite{NO} (see
also \cite{Ols}, \cite{Ner-discrete}), this phenomenon was applied
to construction of exotic unitary representations of $\OO(p,q)$
and $\U(p,q)$.}
It is possible to obtain an  explicit decomposition of the Hilbert space
$\cH_\sigma$ into pieces with uniform spectra, see
\cite{Ner-separation}.

 Next, $\rho_\sigma$ admits a natural analytic continuation
to negative integer $\sigma$,
and the corresponding limit as $\sigma\to-\infty$ is
$ L^2\bigl[\U(p+q)/\U(p)\times\U(q)\bigr]$.
Thus, we obtain some kind of interpolation between
$$
L^2\bigl[\U(p,q)/\U(p)\times\U(q)\bigr]
\qquad\text{and}\qquad
 L^2\bigl[\U(p+q)/\U(p)\times\U(q)\bigr]
 $$

In his short note \cite{Ber2}, 1978, Berezin  gave (without
proofs) the explicit Plancherel
formula
  for representations $\rho_\sigma$
  (in the case of Hermitian symmetric spaces)
  for sufficiently large $\sigma$
  (before the start of bifurcations).
Berezin died soon after this, and the first known proof of his
theorem was  published by Unterberger and Upmeier \cite{UU}, 1994.

In preprints \cite{Ner-beta},  \cite{Ner-berezin}, 1999, the
author concidered the Berezin representations for arbitrary classical
Riemannian symmetric spaces, and obtained the explicit Plancherel
formula   for arbitary $\sigma$.
 \footnote{Some modifications of proofs of  results
 of \cite{Ner-beta} were obtained later in \cite{vDP} and \cite{zha},
see also  \cite{vDH},  \cite{Hil}.}
\footnote{
   Apparently the
Plancherel formula for all the groups can be reduced to
some single
 identity with the
Heckman--Opdam spherical hypergeometric transform (see, for
instance, \cite{HS}), as far as I know, this possibility is  yet
not realized. Apparently also, that the Plancherel formula for Berezin kernels
can be interpolated as an exotic Plancherel formula
 for Heckman--Opdam spherical transform
(this is correct in
rank 1 case, see \cite{Ner-imitation}.}.

Some steps in analysis after the Plancherel formula were undertaked in
\cite{Ner-determinant}, \cite{Ner-shift}, \cite{Ner-imitation}.
Also in \cite{Ner-tits} it was obtained a $p$-adic analog of
Berezin kernels.\footnote{Only for Bruhat--Tits buildings of the
series $A_n$.}

 Since there exists a reach non-$L^2$-analysis
 on Riemannian symmetric space, there arises the following
 question:

 --- {\it are there other homogeneous spaces $G/H$
 and other
kernels $K(x,y)$ admitting an explicit and interesting harmonic analysis?}

Since we discuss the harmonic analysis, a representation $\rho$
having the form (\ref{action}) must be reducible, and hence the distribution
$K(x,y)$ in the formula (\ref{yadro}) is not uniquely  defined. At a
first sight,  a choice is too large.

 Indeed, for each function
 $K(z,u)$ on $\cB_{p,q}\times \cB_{p,q}$
 satisfying
 \begin{equation}
 L\bigl[(a+zc)^{-1} (b+zd),\,(a+uc)^{-1} (b+ud)\bigr]=
L[z,u]
\label{invariance}
,\end{equation}
the expression
\begin{equation}
\langle F_1, F_2\rangle =
 \iint\limits_{\cB_{p,q}\times\cB_{p,q}}
 L(z,u)   F_1(z)\,\ov {F_2(u)}\,d\sigma(z)\,d\sigma(u)
     \end{equation}
     is invariant with respect to the translation
     operators
$$F(z)\mapsto F\bigl( (a+zc)^{-1}(b+zd)\bigr).$$

A kernel $L$ satisfying (\ref{invariance}) can be represented in the form
$$L(z,u)=\ell( \lambda_1,\dots,\lambda_p),$$
where $\lambda_j$ are the eigenvalues of the matrix
$$
(1-zz^*)^{-1} (1-zu^*)(1-uu^*)^{-1}(1-uz^*)
$$
and $\ell$ is a symmetric function on the octant $\lambda_j>0$.
There are many such functions, and hence there are too much
invariant kernels $L$.

Nevertheless, the existence of explicit harmonic analysis is a
non-formal, but very strong restriction; some experience in the
analysis on rank 1 symmetric spaces
(or reading the book \cite{PBM}) shows that a collection  of
possibilities is not  large%
\footnote{There is lot of papers of Molchanov on this subject, for
instance \cite{Mol-comp}, \cite{Mol-sl3}, \cite{vDM}}; for
spaces of rank $>1$ situation became strongly rigid and even an
existence of  examples is not
obvious.

\smallskip

{\bf \punct Natural kernels on pseudo-Riemannian symmetric spaces.}
In the first approximation, the problem that we formulate
is a problem of restriction of  Stein--Sahi representations
to some special symmetric subgroups.

More precisely, for any classical symmetric space
$G/H$ there exists an overgroup $\wt G\supset G$
acting locally on $G/H$ 
(\cite{Ner-uniform}). If $\wt G$ 
is in the list (\ref{line1}--\ref{line3}) (and its action on $G/H$ is locally equivalent
to a fractional--linear  action on a matrix space $\cM$,
see below Section  4),
 then our problem is
a problem of restriction of a Stein representation of $\wt G$ to
the subgroup $G$. In fact, we can forget about the group
$\wt G$ and consider only the restriction $K(\cdot,\cdot)$
 of the Stein kernel to
the symmetric space $G/H$ and action of $G$ in the Hilbert space
defined by this kernel. For some value of the parameter $s$, we obtain
the natural action of $G$ in $L^2(G/H)$.

It is possible to define the kernel $K(\cdot,\cdot)$ in the
terms of the symmetric spaces $G/H$ themselves;
 this allows to extend our problem to the
case, then  $\wt G$ is not in the list (\ref{line1})--(\ref{line3}).
But existence of an island of positivity
of the form
 in these cases becomes a problem.

\smallskip

{\bf \punct Structure of the paper.}
In Section 1, we dicuss our kernels
for the space $G/H\simeq\U(n)\times\U(n)/\U(n)$.
This also gives a realization of the
Stein--Sahi representations of $\U(n,n)$
(For this case, more detailed discussion
including the unipotent representations
   is contained in
\cite{Ner-sobolev}).

In Section 2, we repeat the construction for
spaces $\U(n)/\OO(n)$, $\U(2n)/\Sp(n)$
and for their overgroups
 $\Sp(2n,\R)$, $\SOS(4n)$.

The proof of positivity is contained Section 3,
it is based on the explicit expansions of the kernels in spherical
functions. This done uniformly using  Kadell's
generalized Selberg integral.

In Section 4, we briefly discuss the remaining series
of classical groups. 

In Section 5, we formulate in details
 the problem discussed above in 0.5.
 We also try to explain why the problem
 looks as solvable.

\smallskip

{\bf \punct Acknowledgements.}
I thank V.F.Molchanov, G.I.Olshanski, T.Koba\-yashi, E.van den Ban,
and T.Oshima for discussion of the subject.

This work was done during my visit to
the  Erwin Shr\"odinger Instirute,
Vienna, in winter 2001--2002. I thank
the administrators of the Institute for hospitality.

\medskip

{\bf \punct Notations.}

Let $A=\{a_{ij}\}$ be a matrix. Then $A^t$ is the transposed
matrix, and $A^*$ is the adjoint matrix, i.e., matrix with
elements $\ov a_{ji}$. Also, we denote by $\ov A$ the element-wise
conjugate matrix, i.e.,
the matrix with the matrix elements $\ov
a_{ij}$.

A matrix $x$ over $\R$, $\C$ is {\it symmetric} 
if $x=x^t$ and {\it skew-symmetric} if $x=-x$.

A matrix over $\R$,  $\C$, or quaternions $\H$
is {\it Hermitian} if $x=x*$
and {\it anti-Hermitian} if $x=-x^*$.

If $A$ is a Hermitian matrix, the notation
$A>0$ means that $A$ 
is positive definite, i.e., $hAh^*>0$ for each 
vector-row $h$.

The symbol $\|z\|$ denotes the {\it norm}
of an operator $z$ in a Euclidean space.
The value $\|z\|^2$ coincides with the maximal eigenvalue
of $z^*z$.

\smallskip

We use the standard notation for the Pochhammer symbol
$$(a)_k:=a(a+1)\dots(a+k-1)$$

\smallskip

For a complex number $x$ we use the following notation
for its powers
\begin{equation}
z^{\{\tau\|\sigma\}}:=z^\tau \ov z^\sigma
\end{equation}

\smallskip

The symbol $[a]$ denotes the integer part of $a$.

\bigskip

\begin{center}

{\bf \large 1.  Sobolev kernels on the spaces $\U(n)$

and unitary representations of $\U(n,n)$.}

\end{center}

\medskip

\addtocounter{sec}{1} \setcounter{equation}{0} \setcounter{punct}{0}
\setcounter{theorem}{0}

On construction of this section, see also \cite{Ner-sobolev}.

\smallskip

{\bf \punct  The group $\U(n)$.} The unitary group $\U(n)$ is a
compact Riemannian symmetric space
$$
G/K=\U(n)\simeq \U(n)\times\U(n)/\U(n)
.$$
 Indeed, the group $G:=\U(n)\times\U(n)$
acts on $\U(n)$ by left and right multiplications
$$
(h_1,h_2):\,\,z\mapsto
h_1^{-1} z h_2.
$$
The stabilizer $K\simeq\U(n)$ of the point $z=1$ consists of
elements $(h,h)\in \U(n)\times\U(n)$. Using transformations
$z\mapsto h^{-1}zh$, each element of $\U(n)$ can be reduced to the
diagonal form
 \begin{equation}
\begin{pmatrix}
e^{i\phi_1}&0&\dots&0\\
0&e^{i\phi_2}&\dots&0\\
\vdots&\vdots&\ddots&\vdots\\
0& 0&\dots&e^{i\phi_n}
\end{pmatrix}
\label{can-form-u}
.\end{equation}

\smallskip

{\bf \punct Overgroup.} The space $G/K$ admits a larger group
of symmetries, namely the pseudounitary
 group $\U(n,n)$. Recall that $\U(n,n)$ is the
group of $(n+n)\times (n+n)$ matrices
 $g=\begin{pmatrix} a&b\\c&d\end{pmatrix}$
preserving Hermitian form $H(\cdot,\cdot)$ with the matrix
$\begin{pmatrix} 1&0\\0&-1\end{pmatrix}$, i.e.,
$$H(v\oplus w, v'\oplus w'):=\sum_{j=1}^n v_j\ov v_j'-\sum_{j=1}^n w_j\ov w_j'.$$
 In other
words, the matrix $g$ satisfies the condition
\begin{equation}
\begin{pmatrix} a&b\\c&d\end{pmatrix}^*
\begin{pmatrix} 1&0\\0&-1\end{pmatrix}
\begin{pmatrix} a&b\\c&d\end{pmatrix}
=
\begin{pmatrix} 1&0\\0&-1\end{pmatrix}
\label{def-u-n-n}
. \end{equation}

The group $\U(n,n)$ acts on the space $\U(n)$ by linear fractional
transformations
\begin{equation}
z\mapsto z^{[g]}:=(a+zc)^{-1}(b+zd)
\label{linear-fractional}
.\end{equation}

\begin{lemma}  Let $z\in\U(n)$, $g\in\U(n,n)$.
Then $z^{[g]}\in\U(n)$.
\end{lemma}

{\bf \punct Proof of Lemma 1.1. Identification of $\U(n)$ with a
Grassmannian.}
 For $z\in\U(n)$ consider its graph
 $V_z\subset \C^n\oplus\C^n$. A vector
 $v\oplus w\in\C^n\oplus\C^n $
is an element of $V_z$ if $w=vz$.

For $v,v'\in\C^n$, we have
$$
H(v\oplus vz,v'\oplus v'z)= \langle v,v'\rangle_{\C^n}-\langle
vz,v'z\rangle_{\C^n}= \langle v,v'\rangle_{\C^n}-\langle
v,v'\rangle_{\C^n} =0
,$$
where $\langle \cdot,\cdot\rangle_{\C^n}$ denotes the standard
inner product in $\C^n$.
Hence $V_z$ is a maximal $H$-isotropic subspace%
\footnote{Consider a linear space $L$ equipped with bilinear
(sesquilinear) form $B(\cdot,\cdot)$. A subspace $M\subset L$ is
{\it isotropic} if $B(h,h')=0$ for all $h$, $h'\in M$.} in
$\C^n\oplus \C^n$.

Conversely, let $V$ be a maximal $H$-isotropic subspace in
$\C^n\oplus\C^n$. Since the form $H$ is strictly negative on
$0\oplus \C^n$, we have $V\cap (0\oplus\C^n)=0$. But $\dim V=n$,
and hence $V$ is a graph of some linear operator $z:\C^n\oplus
0\to 0\oplus\C^n$. By the isotropy condition,
$$
0=H(v\oplus vz,v'\oplus v'z)=
H(v\oplus 0,v'\oplus 0)- H(0\oplus vz,0\oplus v'z)
$$
and hence $z\in\U(n)$.

Thus we obtain the natural identification
$$
\left\{
\begin{aligned}
\text{The Grassmannian of maximal}\\
\text{$H$-isotropic subspaces in $\C^n\oplus \C^n$}
\end{aligned}
\right\}
\longleftrightarrow
\left\{
\begin{aligned}
\text{The space $\U(n)=$}
\\
=\U(n)\times\U(n)/\U(n)
\end{aligned}
\right\}
$$

The group $\U(n,n)$ preserves the Hermitian form $H(\cdot,\cdot)$
and hence it transfer isotropic subspaces to isotropic subspaces.
Thus $\U(n,n)$ acts in a natural way on $\U(n)$. It remains to write
an explicit formula for this action.

\smallskip

\begin{lemma}
We have
$g V_{z}= V_{u}$,
where
$
u=z^{[g]}:=(a+zc)^{-1}(b+zd) \label{olsh-notation}
$.
\end{lemma}

 {\sc Proof.} We have $v\oplus vz\in V_z$. Applying $g$, we
obtain
$$v(a+zc)\oplus v(c+zd)\in g V_{z}.$$
Denoting $w=v(a+zc)$, we obtain the required statement.
\hfill $\square$

\smallskip

{\bf \punct Jacobians.}

\begin{lemma}
Denote  by $\mu(z)$ the Haar measure on the group $\U(n)$, denote
by $\mu(z^{[g]})$ its image under the transformation $z\mapsto
z^{[g]}$, given by (\ref{olsh-notation}). Then
\begin{equation}
\mu(z^{[g]})=|\det(a+zc)|^{-2n}\mu(g) \label{jacobian}
\end{equation}
.\end{lemma}

{\sc Proof.} Let $h\in\U(n)\times\U(n)\subset\U(n,n)$ be an
element of the form
\begin{equation} h=\begin{pmatrix}
u&0\\0&v\end{pmatrix} \label{h-in-U}
.\end{equation}
 Then $z^{[h]}=u^{-1}zv$. The Haar measure is
invariant with respect to such transformations.

Denote
\begin{equation}
 \gamma(g,z):=|\det(a+zc)|^{-2n}
 \label{multip}
  .\end{equation}
 This expression satisfies the
identity
$$
\gamma(h_1 g h_2, z^{[h_2^{-1}]})=\gamma (g,z)
$$
for $h_1$, $h_2$ having the form (\ref{h-in-U}). Due the
invariance of the Haar measure, the Jacobian (or Radon-Nykodim
derivative)
$$
\mu(z^{[g]})/\mu(z)
$$
satisfies the same condition.

Hence, without loss of generality, we can assume that $z=1$,
$z^{[g]}=1$. The tangent space to $\U(n)$ at this point consists
of matrices $\Delta$ satisfying $\Delta+\Delta^*=1$.
We denote this space of all the anti-Hermitian
matrices by $\AHerm$.

It is easy to show (see, for instance, \cite{Ner-beta}, Lemma 1.1),
that the differential of a linear-fractional transformation
$z\mapsto z^{[g]}$ is given by
\begin{equation}
\Delta\mapsto (a+zc)^{-1} \Delta (-c z^{[g]}+d)
\label{differential}
,\end{equation}
this formula is valid on whole space of $n\times n$ matrices for
$g\in\GL(2n,\C)$, in particular it can be used in our case.


Now $z=z^{[g]}=1$ and hence $a+c=b+d$.
This equation, together with (\ref{def-u-n-n})
allows to transform 
(\ref{differential}) to the form

\begin{equation}
\Delta\mapsto  (a+c)^{-1} \Delta  \bigl[(a+c)^{-1}\bigr]^*
\label{1}
.\end{equation}


It remains to evaluate  the determinant $\delta(P)$ of the linear
transformation $\AHerm(n)\to\AHerm(n)$ given by
$$
\Delta \mapsto P\Delta P^*, \qquad P\in \GL(n,\C)
. $$
 Obviously, $\delta(P)$ is a homomorphism
from $\GL(n,\C)$ to the multiplicative group of positive numbers.
Hence $\delta(P)=1$ on $\U(n)$. Each $P\in \GL(n,\C)$ can be
represented in the form $u\Lambda v$, where $u$, $v\in\U(n)$, and
$\Lambda$ is a diagonal matrix. Obviously
$\delta(\Lambda)=|\det(\Lambda)|^{2n}$ and hence
$\delta(P)=|\det(P)|^{2n}$. This finishes the proof. \hfill
$\square$

\smallskip

{\bf \punct Invariant kernels on $\U(n)$.} Denote by $\cB_n$ the
set of $n\times n$ matrices $z$ with norm $<1$; an equivalent
condition is $1-z^*z>0$
 By $\ov \cB_n$ denote the set of all
$n\times n$ matrices with norm $\le 1$; an equivalent condition is
$1-z^*z\ge 0$.

Fix $\sigma\in\C$. For $z\in\cB_n$, we define
\begin{equation}
(1-z)^\sigma:= \sum_{j=0}^\infty \frac{(-\sigma)_j}{j!} z^j
 \label{power}
,\end{equation}
 the series in the right-hand side is convergent.

Fix $\sigma$, $\tau\in\C$. We define the function
$\det(1-z)^{\{\sigma|\tau\}}$ depending in the variable
$z\in\cB_n$ by
\begin{equation}
\det(1-z)^{\{\sigma|\tau\}} = \det\Bigl[(1-z)^\sigma\Bigr]
\det\Bigl[(1-\ov z)^\tau\Bigr]
 \label{double-power}
.\end{equation}
This expression is a smooth $\C$-valued function on $\cB_n$.
 For
$z\in\ov\cB_n$, we define
\begin{equation}
\det(1-z)^{\{\sigma|\tau\}}:= \lim_{\epsilon\to +0}
\det(1-(1-\epsilon)z)^{\{\sigma|\tau\}}
\label{23}
.\end{equation}
 This limit exists outside the surface $\det(1-z)=0$,
and the expression (\ref{23}) is a smooth function on the set
$$
1-zz^*\ge 0,\qquad \det(1-z)\ne 0
.$$

Obviously, for a unitary matrix $h\in\U(n)$,
$$
\det(1-hzh^{-1})^{\{\sigma|\tau\}}=\det(1-z)^{\{\sigma|\tau\}}
.$$

\begin{lemma}
For $z\in\U(n)$, denote by $e^{i\psi_1}$,\dots, $e^{i\psi_n}$ its
eigenvalues. Then
$$
\det(1-z)^{\{\sigma|\tau\}}=
 \exp\Bigl\{(\sigma-\tau)\sum_{k=1}^n(\psi_k-\pi)/2\Bigr\}
\cdot\prod_{k=1}^n |\sin (\psi_k/2)|^{\sigma+\tau}
.$$
\end{lemma}

This statement is more-or-less obvious, 
for a formal proof, see \cite{Ner-sobolev}

\smallskip

We define the kernel $L_{\sigma,\tau}(z,u)$ on $\U(n)$ by
\begin{equation}
L_{\sigma,\tau}(z,u)=\det(1-zu^*)^{\{\sigma|\tau\}}, \qquad
z,u\in\U(n)
 \label{kernel-u}
.\end{equation}

\begin{lemma}

{\rm a)}  For $\tau$, $\sigma\in\R$, the kernel
$L_{\sigma,\tau}(z,u)$ is Hermitian, i.e.,
$$
L_{\sigma,\tau}(u,z)=\ov{L_{\sigma,\tau}(z,u)}
.$$

{\rm b)} The kernel $L_{\sigma,\tau}(z,u)$
 is $\U(n)$-invariant, i.e.,
 $$
L_{\sigma,\tau}(h_1zh_2,h_1uh_2)=L_{\sigma,\tau}(z,u), \qquad
h_1,h_2\in\U(n)
$$
\end{lemma}

The both statements are obvious.

\smallskip

 We  define the sesquilinear form $I_{\sigma,\tau}$ on
$C^\infty((\U(n))$ given by
\begin{equation}
I_{\sigma,\tau}(F_1,F_2)=
\iint_{\U(n)\times\U(n)}
 L_{\sigma,\tau}(z,u)
F_1(z) \ov{ F_2(u)}\,d\mu(z)\,d\mu(u) \label{form-u}
\end{equation}

If $\sigma$, $\tau\in \R$, then the form is Hermitian.
 The integral
is convergent if $\Re(\sigma+\tau)>-1$.
 By general reasons, the
integral admits a meromorphic continuation to the whole plane
$(\sigma,\tau)\in\C^2$ (see \cite{Bern}, \cite{Ati}). In fact, our
proof of positivity given below is based on
an the explicit construction
of this meromorphic continuation.

\smallskip

{\bf \punct Action of $\U(n,n)$.}

\nopagebreak

\begin{lemma}
The operators $\rho_{\sigma,\tau}$ given by
\begin{equation}
\rho_{\sigma,\tau}\begin{pmatrix} a&b\\c&d\end{pmatrix} F(z)=
F\bigl[(a+zc)^{-1}(b+zd)\bigr] \det(a+zc)^{\{-n-\sigma|-n-\tau\}}
\label{representation-u}
\end{equation}
  preserve
the form $I_{\sigma,\tau}$
\end{lemma}

Before starting a proof, we  give some preliminary comments on
formula (\ref{representation-u}).

\smallskip

{\bf 1.} We must define these operators more carefully. First,
$$\det(a+zc)^{\{-n\|-n\}}=|\det(a+zc)|^{-2n}$$
and this expression is well-defined. Further
\begin{equation}
\det(a+zc)^{\{-\sigma\|-\tau\}}= \det(a)^{\{-\sigma\|-\tau\}}
\det(1+zca^{-1})^{\{-\sigma\|-\tau\}} \label{2}
\end{equation}
Since $\|ca^{-1}\|<1$ (this easily follows from (\ref{def-u-n-n})),
we have also $\|zca^{-1}\|<1$. Thus the factor
$\det(1+zca^{-1})^{\{-\sigma\|-\tau\}}$ is well defined, see
(\ref{power}), (\ref{double-power}). It remains to define
$$
\det(a)^{\{-\sigma\|-\tau\}}= |\det(a)|^{-\sigma-\tau}
\exp\bigl\{i(\tau-\sigma)\arg (a) + 2\pi i(\tau-\sigma)k\bigr\}
,$$
where $k$ ranges in $\Z$. If $(\sigma-\tau)\in \Z$, then the last
expression is well defined, and $\rho_{\sigma,\tau}$ is a linear
representation of $\U(p,q)$.

Otherwise, the multi-valued function (\ref{2}) on $\cB_n$
 splits into a
countable family of  smooth branches%
\footnote{In particular, our operators preserve the space
$C^\infty(\U(n)$.)}; these branches differs by constant factors
$\exp\{2\pi i(\tau-\sigma)k\}$. Thus the formula
(\ref{representation-u}) for a given $g$,
  defines a
countable family of operators $\rho_{\sigma,\tau}$, which differ one from another by
constant factors. We can choose one such operator in an arbitrary
way, and then we will obtain a unitary {\it projective} (see, for
instance, \cite{Kir}, Section 14) representation of $\U(n,n)$,
\begin{equation}
\rho_{\sigma,\tau}(g)\rho_{\sigma,\tau}(h)=
\kappa(g,h)\rho_{\sigma,\tau}(gh),\qquad\kappa(g,h)\in\C^*
\label{projective-representation}
.\end{equation}

{\bf 2.} Since $\det(a)$ does not vanish, the function
$\ln\det(a)$ is a well-defined function on the universal covering
$\U(n,n)^\sim$ of $\U(n,n)$. Hence the expression (\ref{2}) is a
well defined single-valued expression on $\U(n,n)^\sim$; so we  can
consider $\rho_{\sigma,\tau}$ as a linear representation of
$\U(n,n)^\sim$,
\begin{equation}
\rho_{\sigma,\tau}(g)\rho_{\sigma,\tau}(h)=
\rho_{\sigma,\tau}(gh)
 \label{linear-representation}
.\end{equation}

{\bf 3.} We must prove the identities
(\ref{projective-representation}), (\ref{linear-representation}).
A direct calculation is not difficult, but it is more reasonable
to avoid it. Denote
$$\nu(g,z):=\det(a+zc)^{\{-n-\sigma|-n-\tau\}}.$$
The desired identities are reduced to the "cocycle identity"
\begin{equation}
\nu(g h, z)=\nu(h,z)\nu(g,z^{[h]})
\label{cocycle-identity}
.\end{equation}
 But $\nu(g,h)$ is a power of the
expression $\gamma(g,z)$ given by (\ref{multip}), the latter
expression is the Jacobian. For a Jacobian, the cocycle identity
is obvious.

\smallskip

{\bf 4.} Let $\tau=0$. Then our construction gives a
representation of the Harish-Chandra holomorphic  series.
The kernel $L_{\sigma,0}$ is the standard reproducing kernel for
these representations, see, for instance,
 \cite{Ber}, \cite{VR}, \cite{NO}. 
Representations of holomorphic series
 admit realizations (see \cite{NO})
in spaces of
holomorphic functions, in spaces of distributions on matrix balls
$\cB_n$, and in spaces of distributions on Shilov boundary $\U(n)$
of the matrix ball. The last variant corresponds to our
realization.

For $\sigma=0$, we obtain a lower weight representation.

\smallskip

{\bf 5.} Our construction is invariant with respect to the shift
$$(\sigma,\tau)\mapsto(\sigma+1,\tau-1).$$
We will not refer to this remark below, for discussion, see
\cite{Ner-sobolev}, 2.9.

\smallskip

  {\bf \punct Proof of Lemma 1.6.} We use the following simple identity
$$
\det\bigl(1-z^{[g]}\bigl[u^{[g]}\bigr]^*\bigr) =
\det(1-zu^*)\det(a+zc)^{-1}
\ov{\det(a+uc)}^{-1}
$$
valid for $g=\begin{pmatrix}a&b\\c&d\end{pmatrix}\in\U(n,n)$.
 Keeping it in mind and using the formula
 (\ref{jacobian}) for the Jacobian,
 we substitute $z\mapsto z^{[g]}$, $u\mapsto u^{[g]}$
to integral (\ref{form-u}) and obtain
\begin{multline*}
\iint_{\U(n)\times\U(n)}
 L_{\sigma,\tau}(z,u)\det(a+zc)^{\{-\sigma|-\tau\}}
 \ov{\det(a+uc)}^{\{-\sigma|-\tau\}}
\times\\ \times F_1(z^{[g]}) \ov{
F_2(u^{[g]})}\det(a+zc)^{-2n}\det(a+uc)^{-2n}\,d\mu(z)\,d\mu(u)
\end{multline*}
Q.E.D.

  {\bf \punct Unitary representations of $\U(n,n)$.}

\begin{theorem}
Let $\sigma$, $\tau\notin \Z$. The Hermitian form
$I_{\sigma,\tau}$ given by (\ref{form-u}) is definite
iff\footnote{recall that $[x]$ denotes the integral part of $x$.}
$$[-\tau]=[\sigma+n].$$
\end{theorem}

Proof is given below in Section 3. More elementary proof is
contained in \cite{Ner-sobolev}. This paper also contains a
picture of a domain of positivity.

\smallskip

It is convenient to introduce new parameters $t$, $s$ by
\begin{equation}
\sigma=-n/2+s,\qquad \tau=-n/2+t \label{ts-u}
.\end{equation}

In this notation, the cases of even and odd $n$ are slightly
different.

a) For an odd $n$ the condition of positivity is
$$
|s-j|<1/2,\qquad |t-j|<1/2,\qquad
\text{for some $j\in \Z$}
.$$

b) For an even $j$ the condition of positivity is
$$t\in [j,j+1],\qquad s\in [j-1,j]
\qquad \text{for some $j\in \Z$}
.$$

\smallskip

\begin{prop}
If $s=t$, then the Hermitian form $I_{\sigma,\tau}$ coincides with
the $L^2$-inner product
$$\langle F_1,F_\rangle
=\int_{\U(n)} F_1(z)\ov{F_1(z)}d\mu(z)
.$$
\end{prop}

See a proof in \cite{Ner-sobolev}, this also can be easily derived
from calculations of our Section 3.

Under the conditions of Theorem 1.7, we denote by
$\cH_{\sigma,\tau}$ the completion of $C^\infty(\U(n))$ with
respect to this form. We obtain that our representation
$\rho_{\sigma,\tau}$ in this case is unitary in $H_{\sigma,\tau}$.

\bigskip

{\bf \large 2.  Sobolev kernels on the spaces
 $\U(n)/\OO(n)$,
$\U(2n)/\Sp(n)$ and unitary representations of the groups
$\Sp(2n,\R)$, $\SOS(4n)$}

\medskip

\addtocounter{sec}{1} \setcounter{equation}{0}
\setcounter{punct}{0} \setcounter{theorem}{0}

Here we show that the space $\U(n)/\OO(n)$ can be realized as the
space of unitary  symmetric  matrices and the
space $\U(2n)/\Sp(n)$ as the space of unitary skew-symmetric
$2n\times 2n$-matrices %
\footnote{A remark for experts. These spaces are realized as the
Shilov boundaries of the bounded Cartan domains
$\Sp(2n,\R)/\U(n)$, $\SO^*(4n)/\U(2n)$, see \cite{PS}.}.
 We show
that the groups $\Sp(2n,\R)$ and $\SOS(4n)$ respectively act on
these spaces by linear-fractional transformations. Then we
restrict the kernel $L_{\sigma,\tau}$ defined by (\ref{kernel-u})
and obtain unitary representations of $\Sp(2n,\R)$ and
$\SO^*(4n)$.

\bigskip

{\bf A. Spaces $\U(n)/\OO(n)$.}

\medskip

{\bf \punct Symmetric spaces $\U(n)/\OO(n)$.} Now let $\cX_n$ be
the space $n\times n$ matrices $z$ satisfying the conditions
$$zz^*=1,\qquad z=z^t.$$

The group $\U(n)$ acts on $\cX_n$ by the transformations
$$
h:\,\,z\mapsto h^t z h
.$$
The stabilizer of the point $z=1$  is the real orthogonal group
$\OO(n)$. It can be easily verified that action of $\U(n)$ on
$\cX_n$ is transitive. Thus,
$$\cX_n\simeq \U(n)/\OO(n).$$

Each element of $\cX_n$ can be reduced by the transformations
$z\mapsto h^tzh$, where $h\in \OO(n)$, to the form
\begin{equation}
\begin{pmatrix}
e^{i\phi_1}&0&\dots&0\\
0&e^{i\phi_2}&\dots&0\\
\vdots&\vdots&\ddots&\vdots\\
0& 0&\dots&e^{i\phi_n}
\end{pmatrix}
\label{can-form-uo}
.\end{equation}
 The collection $\phi_j$ is uniquely
defined up to permutations.

A matrix $\Delta$ is an element of the tangent space to $\cX_n$ at
$z=1$ iff
$$\Delta=\Delta^t,\qquad\text{and $\tfrac 1i \Delta$ is a real
matrix}
.$$

{\bf\punct Overgroup.} Again, our space $\cX_n$ admits a larger
group of symmetries, namely $\Sp(2n,\R)$.

To observe this, we realize the real symplectic group $\Sp(2n,\R)$
as a group of {\it complex} $(n+n)\times(n+n)$ matrices $g$
preserving the Hermitian form in $\C^n\oplus\C^n$ with the matrix
$\begin{pmatrix}1&0\\0&-1\end{pmatrix}$ and the skew-symmetric
bilinear form $B$ having the matrix
$\begin{pmatrix}0&1\\-1&0\end{pmatrix}$. In other words,
$$
g^*\begin{pmatrix}1&0\\0&-1\end{pmatrix}g=
\begin{pmatrix}1&0\\0&-1\end{pmatrix};
\qquad g^t\begin{pmatrix}0&1\\-1&0\end{pmatrix}g=
\begin{pmatrix}0&1\\-1&0\end{pmatrix}
.$$
These two identities also imply that the matrix $g$ has the
following block structure
$$
g=\begin{pmatrix}\Phi&\Psi\\ \ov \Psi&\ov \Phi\end{pmatrix}
.$$
Equivalently $g$ preserves the  real subspace of $\C^n\oplus\C^n$
consisting of vectors $h\oplus \ov h$.

\smallskip

The group $\Sp(2n,\R)$ acts on $\cX_n$ by linear fractional
transformations
$$
g: \,\, z\mapsto z^{[g]}:=(\Phi+z\ov \Psi)^{-1}(\Psi+z\ov\Phi)
.$$

\begin{lemma}
 Let $z\in\cX_n$, $g\in\Sp(2n,\R)$.
Then $z^{[g]}\in \cX_n$.
\end{lemma}

\smallskip

{\sc Proof.} As above, for a matrix $z\in\cX_n$, we consider its
graph $V_z\subset \C^n\oplus\C^n$. Literally repeating the
considerations of Subsection 1.2, we obtain that the subspace
$V_z$ is maximal isotropic with respect to the Hermitian form $H$.

The condition $z=z^t$ means that $V_z$ is Lagrangian (=maximal
isotropic) with respect to the skew-symmetric form $B$.
Indeed,
$$
B(v\oplus vz,w\oplus wz)=vzw-wzv
$$
where $v$, $w$ are matrices-rows. Since $z=z^t$,
we obtain zero in the right-hand side.

Thus, we obtain the identification
\begin{multline*}
\left\{
\begin{aligned}
\text{Grassmannian of $n$-dimensional subspaces
 in $\C^n\oplus\C^n$}
 \\
 \text{isotropic with respect
  to the both forms $H$ and $B$}
 \end{aligned}
 \right\}
 \longleftrightarrow\\
\longleftrightarrow
 \left\{
 \begin{aligned}
\text{The space $\cX_n=\U(n)/\OO(n)$ of $n\times n$ matrices}
\\
\text{$z$ satisfying $z=z^t$, $z^*z=1$}
\end{aligned}
\right\}
.\end{multline*}

An element $g\in\Sp(2n,\R)$ preserves the both forms $H$, $B$ and
hence it maps our isotropic Grassmannian to itself. \hfill
$\square$.

We denote by $\mu$ the unique $\U(n)$-invariant measure on
$\cX_n$.
 The
Jacobian of the transformation $z\mapsto z^{[g]}$ is given by the
following formula
\begin{equation}
\mu(z^{[g]})=\det(\Phi+z\Psi)^{-(n+1)} \mu(z) \label{jacobian-uo}
\end{equation}

A proof is the same as above (Lemma 1.3).

\smallskip

{\bf \punct Invariant kernels on $\U(n)/\OO(n)$.} Now fix $\sigma$,
$\tau\in\C$. Consider the kernel $L_{\sigma,\tau}(z,u)$, where
$z$, $u\in\cX_n$ 
the same kernel (\ref{kernel-u})
 restricted to $\cX_n$.

Obviously, the kernel $L_{\sigma,\tau}(z,u)$
 is $\U(n)$-invariant
 $$
 L_{\sigma,\tau}(h^tzh,h^tuh)=L_{\sigma,\tau}(z,u)
 $$
 We define the sesquilinear form
 on $C^\infty(\cX_n)$ by
\begin{equation}
I_{\sigma,\tau}(F_1,F_2)= \iint_{\cX_n\times\cX_n}
 L_{\sigma,\tau}(z,u)
F_1(z) \ov{ F_2(u)}\,d\mu(z)\,d\mu(u) \label{form-uo}
.\end{equation}

 \smallskip

Again, we consider the meromorphic continuation of
$I_{\sigma,\tau}$ to the domain $(\sigma,\tau)\in\C^2$.

As above, we define the linear operators
$$
\rho_{\tau,\sigma}
\begin{pmatrix}\Phi&\Psi\\ \ov \Psi&\ov \Phi\end{pmatrix}
f(z)= f\bigl( (\Phi+z\ov\Psi)^{-1} (\Psi+z\ov\Phi)\bigr)
\,\det(\Phi+z\ov\Psi)^{\{-(n+1)/2-\sigma\|-(n+1)/2-\tau\}}
.$$
As above, these operators preserve the form $I_{\sigma,\tau}$.

\begin{theorem}
Let $n>1$. Let $2\sigma$, $2\tau\notin \Z$. The Hermitian form
$I_{\sigma,\tau}$ is definite iff
$$[-2\tau-n-1]=[2\sigma].$$
%
%
%
\end{theorem}

Proof is given below in Section 3.

It is convenient to define new parameters $s$, $t$ by
$$\sigma=-(n+1)/4+s,\qquad \tau=-(n+1)/4+t.$$

The conditions of positivity are

a) For an even $n$,
$$s,t\in[j/2-1/4,j/2+1/4]\qquad\text{for some $j\in\Z$}
.$$

b) For an odd $n$
$$s\in [j/2,j/2+1/4],\qquad
t\in [j/2-1/4, j/2] \qquad\text{for some $j\in\Z$}
.$$

For $s=t$ the form $I_{\sigma,\tau}$
 defines the $L^2$-inner product.

\bigskip

{\bf B. Spaces $\U(2n)/\Sp(n)$.}

\medskip

{\bf \punct The space $\U(2n)/\Sp(n)$.} Now we consider the space
$\cY_n$ of $2n\times 2n$ matrices $z$ satisfying the conditions
$$zz^*=1,\qquad z=-z^t.$$
An example of $z$ satisfying these conditions is
\begin{equation}
J=\begin{pmatrix}0&1\\-1&0\end{pmatrix} \qquad\text{where $1$ is
the $n\times n$ unit matrix}
 \label{matrix-J}
.\end{equation}

  The group $\U(2n)$ acts on this
space $Y_n$ by the transformations
$$h:\qquad z\mapsto h^tzh.$$

The stabilizer of the point $J\in\cY_n$ consists of matrices
satisfying
$$hh^*=1,\qquad h^t Jh=J.$$
These  equations give one of the
 standard realizations of the compact
symplectic group $\Sp(n)$. It is easy to show that $\cY_n$ is a
homogeneous space
$$\cY_n=\U(2n)/\Sp(n).$$
Using the transformations $z\mapsto  h^tzh$, where $h\in\Sp(n)$,
we can reduce each element of $\cY_n$ to the form
\begin{equation}
\Omega=\begin{pmatrix} 0&\Lambda\\-\Lambda&0
\end{pmatrix},\qquad
\text{where}\qquad \Lambda=
\begin{pmatrix}
e^{i\phi_1/2}&0&\dots&0\\
0&e^{i\phi_2/2}&\dots&0\\
\vdots&\vdots&\ddots&\vdots\\
0& 0&\dots&e^{i\phi_n/2}
\end{pmatrix}
\label{can-form-usp}
.\end{equation}
 The numbers are uniquely defined up to
permutations and transformations $\phi_j\mapsto-\phi_j$.

\smallskip

{\bf\punct Overgroup.}
We realize the  group $\SOS(4n)$
 as a group of {\it complex} $(2n+2n)\times(2n+2n)$ matrices $g$
preserving the Hermitian form in $\C^{2n}\oplus\C^{2n}$ with the
matrix $\begin{pmatrix}1&0\\0&-1\end{pmatrix}$
 and the
symmetric bilinear form $B$ having the matrix
$\begin{pmatrix}0&1\\-1&0\end{pmatrix}$. In other words,
$$
g^*\begin{pmatrix}1&0\\0&-1\end{pmatrix}g=
\begin{pmatrix}1&0\\0&-1\end{pmatrix};
\qquad g^t\begin{pmatrix}0&1\\1&0\end{pmatrix}g=
\begin{pmatrix}0&1\\1&0\end{pmatrix}
$$
These two identities also imply that the matrix $g$ has the
following block structure
$$
g=\begin{pmatrix}\Phi&\Psi\\- \ov \Psi&\ov \Phi\end{pmatrix}
$$

The group $\SOS(4n)$ acts on $\cY_n$ by 
the linear fractional
transformations
\begin{equation}
g: \,\, z\mapsto z^{[g]}:=(\Phi-z\ov \Psi)^{-1}(\Psi+z\ov\Phi)
\label{linear-fractional-usp}
.\end{equation}

\begin{lemma}
 Let $z\in\cY_n$, $g\in\SOS(4n)$.
Then $z^{[g]}\in \cY_n$.
.\end{lemma}

\smallskip

{\sc Proof} is the same as above.  For a matrix $z\in\cY_n$, we
consider its graph $V_z\subset \C^{2n}\oplus\C^{2n}$. As it was
shown in 1.3, the subspace $V_z$ is maximal isotropic with respect
to the Hermitian form $H$. The condition $z=-z^t$ means that $V_z$
is maximal isotropic with respect to the symmetric form $B$.

Thus, we obtain the identification
\begin{multline*}
\left\{
\begin{aligned}
\text{Grassmannian of $n$-dimensional subspaces
 in $\C^{2n}\oplus\C^{2n}$}
 \\
 \text{isotropic with respect
  to the both forms $H$ and $B$}
 \end{aligned}
 \right\}
 \longleftrightarrow\\
\longleftrightarrow
 \left\{
 \begin{aligned}
\text{The space $\cY_n=\U(2n)/\Sp(n)$ of $2n\times 2n$ matrices}
\\
\text{$z$ satisfying $z=-z^t$, $z^*z=1$}
\end{aligned}
\right\}
\end{multline*}
Now the statement became obvious.

\smallskip

{\bf\punct Jacobian.}

\begin{lemma}
For the linear-fractional transformations
(\ref{linear-fractional-usp}), we have
\begin{equation}
\mu(z^{[g]})=|\det(\Phi-z\ov\Psi)|^{-(4n-2)}\mu(z)
\label{jacobian-usp}
.\end{equation}
\end{lemma}

{\sc Proof.} The both expressions,
\begin{equation}
\nu_1(g,z):=\det(\Phi-z\ov\Psi);\qquad
\nu_2(g,z):=\frac{\mu(z^{[g]})}{\mu(z)}
 \label{two-cocycles}
.\end{equation}
 satisfy the cocycle identity
(\ref{cocycle-identity}) mentioned above.

Recall some standard facts about the solutions of the cocycle
identity, see \cite{Kir}, 13.2. Let $G/R$ be a homogeneous space,
let $u$ be a point fixed by the subgroup $R$. Let $\nu$ be a
solution of the cocycle equation on $G/R$. Obviously, for $r_1$,
$r_2\in\R$,
$$\nu(r_1r_2,u)=\nu(r_1,u)\nu(r_2,u).$$
i.e., the function $r\mapsto\nu(r,u)$ is a character of $R$.
Moreover, this correspondence between solutions of
(\ref{cocycle-identity}) and characters of $R$ is a bijection.

In our case, $G=\SOS(4n)$ and $R$ is a maximal parabolic in $G$,
the reductive part of $R$ is the group $\GL(n,\H)$. Each character
of $\GL(n,\H)$ has the form $r\mapsto \det(r)^s$, and hence our
question is reduced to an evaluation of $s$.

Since we must know the exponent $s$, it is sufficient to find
(\ref{two-cocycles}) only for $z=J$ and some appropriate $g$ lying
in the stabilizer of $J$.

We choose $g=r(t)$ being the block $(n+n+n+n)\times
(n+n+n+n)$-matrix
$$
r(t):=\begin{pmatrix}
 \cosh t &0 &0 & \sinh t\\
 0&\cosh t& -\sinh t&0\\
 0&-\sinh t& \cosh t&0\\
 \sinh t & 0&0&\cosh t
\end{pmatrix}
.$$
A direct calculation shows that $J^{[r(t)]}=J$. Also
\begin{equation}
\det(\Phi-J\Psi)^{-1}=e^{2nt}
\label{4}
. \end{equation}
  Applying
the formula (\ref{differential}), we obtain that  the differential
of $z\mapsto z^{[r(t)]}$ at $J$ is
\begin{equation}
\Delta\mapsto e^{4t}\Delta
 \label{3}
 \end{equation}
 Since $\dim
\cY_n=\dim\U(2n)-\dim\Sp(n)=2n^2-n$, the determinant of the
transformation (\ref{3}) is $\exp(4n(2n-1)t)$. Comparing with
(\ref{4}), we obtain the exponent in (\ref{jacobian-usp}). \hfill
$\square$

\smallskip

{\bf\punct Unitary representations of $\SOS(2n)$.} We define the
canonical invariant kernel on $\cY_n$ as the restriction of the
kernel
$$ L_{\sigma,\tau}(z,u)=\det(1-zu^*)^{\{\sigma\|\tau\}}
$$ to $\cY_n$.
We also define the Hermitian form on $C^\infty(\cY_n))$ by
$$
I_{\sigma,\tau}(F_1,F_2)= \iint_{\cY_n\times\cY_n}
 L_{\sigma,\tau}(z,u)
F_1(z) \ov{ F_2(u)}\,d\mu(z)\,d\mu(u)
.$$

\begin{theorem}
Let $n>1$. Let $\sigma$, $\tau\notin\Z$. The form
$I_{\sigma,\tau}$ is sign definite iff
$$[-\tau]=[\sigma+2n-1].$$
\end{theorem}

The proof is contained in the next section.

\smallskip

We define the representations $\rho_{\sigma,\tau}$ of $\SOS(4n)$
in $C^\infty(\cY_n)$ by
$$\rho_{\sigma,\tau}
\begin{pmatrix}\Phi&\Psi\\-\ov\Psi&\ov\Phi\end{pmatrix}
 F(z)=
F\bigl((\Phi-z\ov\Psi)^{-1}(\Psi+z\ov\Phi)\bigr)
\det(\Phi-z\ov\Psi)^{\{-n+1/2-\sigma\|-n+1/2-\tau\}}
.$$
These operators preserve the form $I_{\sigma,\tau}$.
 Under the
conditions of Theorem 2.5, our representations are unitary.

It is also convenient to introduce new parameters $t$, $s$ by
$$\sigma=-n+1/2+s,\qquad \tau=-n+1/2+t.$$
Then the condition of positivity is
$$
|s-j|<1/2,\qquad |t-j|<1/2,\qquad \text{for some $j\in\Z$}
.$$

\bigskip

{\bf \large 3.  Application of the Kadell integral}

\medskip

\addtocounter{sec}{1} \setcounter{equation}{0}
\setcounter{punct}{0} \setcounter{theorem}{0}

{\bf \punct Jack polynomials. Preliminaries,} for details see
Macdonald, \cite{Mac}, VI.4, VI.10, and further references in this
book. We consider $n$-dimensional torus $\T^n$ as a subset in
$\C^n$ consisting of points
$$(x_1,\dots,x_n),\qquad |x_1|=1,\dots,|x_n|=1.$$
It is also convenient to write
$$x_1=e^{i\phi_1},\dots,x_n=e^{i\phi_n}.$$
We denote by $\Pol_n$ the space of all polynomials in
$x_1$,\dots,$x_n$ symmetric with respect to permutations of $x_j$.
We denote by $\LPol_n$ the space of all symmetric Laurent
polynomials in $x_1^{\pm 1}$,\dots, $x_n^{\pm 1}$.

Fix $\kappa>0$.
 Consider the inner
products in the spaces $\Pol_n$ and $\LPol_n$ given by
\begin{equation}
\langle f,g\rangle_\kappa= \int_{\T^n} f(x)\ov{g(x)} \prod_{1\le
k<l\le n}|x_k-x_l|^{2\kappa} \prod_{k=1}^n d\phi_k
\label{jack-polynomials}
 .\end{equation}
 The Jack polynomials are
orthogonal polynomials with respect to this scalar product.  In a
multivariate case, the Gramm--Schmidt orthogonalization of
polynomials is not a canonical operation, and hence  we must
define them more carefully.

Let $\lambda$ be a collection of integers (a Young diagram)
\begin{equation}
\lambda:\quad \lambda_1\ge\lambda_2\ge\dots\ge\lambda_n\ge 0
\label{young}
.\end{equation}
 Denote
$$|\lambda|=\sum \lambda_j.$$
Let $|\lambda|=|\mu|$. We say $\lambda\ge\mu$ iff
$$
\lambda_1+\dots+\lambda_j\ge\mu_1+\dots+\mu_j
 \qquad \text{for
all $j$}
.$$

We define the monomial symmetric function $m_\lambda$ as
$$
m_\lambda=x_1^{\lambda_1}x_2^{\lambda_2}\dots
x_n^{\lambda_n}+\dots
,$$
where $\dots$ is the sum of all pairwise distinct monomials
$x_{\sigma(1)}^{\lambda_1}\dots x_{\sigma(n)}^{\lambda_n}$, where
$\sigma$ is a permutation.

The {\it Jack polynomials} $P_\lambda(x)\in \Pol_n$ are defined by
two conditions:

a) $P_\lambda$ are orthogonal with respect to the scalar product
(\ref{jack-polynomials})

\begin{equation}
\!\!\!\!\!\!\!\!\!\!\!\!\!\!\!\!\!\!\!\!\!\!\!\!
\!\!\!\!\!\!\!\!\!\!\!\!\!\!\!\!\!\!\!\!\!\!\!\!
\!\!\!\!\!\!\!\!\!\!\!\!\!\!\!\!\!\!\!\!\!\!\!\!\!\!\!
\!\!\!\!\!\!\!\!\!\!\!\!\!\!\!\!\!\!\!\!
 b) \qquad
P_\lambda=m_\lambda
+\sum\limits_{\mu:\,|\mu|=|\lambda|,\mu<\lambda}u_\mu m_\mu
\label{u}
 \end{equation}

Uniqueness of the Jack polynomials is obvious but existence is a
theorem.

\smallskip

The Jack polynomials satisfy the following identity
$$
P_{\lambda_1+1,\dots,\lambda_n+1}(x)= x_1\dots x_n
P_{\lambda_1,\dots,\lambda_n}(x)
.$$
This allows to define Laurent Jack polynomials
$P_\lambda\in\LPol_n$ for
\begin{equation}
\lambda:\quad \lambda_1\ge\lambda_2\ge\dots\ge\lambda_n,
\qquad\lambda_j\in \Z
\label{signatures}
\end{equation}
 (we permit negative $\lambda_j$) by
$$
P_{\lambda,\dots,\lambda_n}(x) =(x_1\dots x_n)^{-m}
P_{\lambda_1+m,\dots,\lambda_n+m}(x)
,$$
where $m$ is sufficiently large ($\lambda_n+m\ge 0$).

This system of polynomials forms an orthogonal basis in $\LPol_n$
with respect to the scalar product (\ref{jack-polynomials}).

\smallskip

{\bf\punct Radial functions.} Consider one of our symmetric spaces
\begin{equation}
K/H=\U(n)\times \U(n)/\U(n),\quad \U(n)/\OO(n), \quad
\U(2n)/\Sp(n) \label{spaces-a}
.\end{equation}
Consider an $H$-invariant function $F$ on
$G/H$. Obviously, this function can be considered as a symmetric
function depending in the invariants $\phi_1$,\dots, $\phi_n$
defined respectively in (\ref{can-form-u}), (\ref{can-form-uo}),
(\ref{can-form-usp}).

Consider the map
 \begin{equation}
z\mapsto (\phi_1,\dots,\phi_n) \label{radial-map} \end{equation}
that takes a matrix $z$ to its collection
of invariants. Since $\phi_j$ are defined up to a permutation, we
must assume
$$
0\le\phi_1\le\phi_2\dots\le \phi_n<2\pi
.$$

The following variant of the Weyl integration formula holds.

\begin{lemma}
The push-forward%
\footnote{Let $A$ be a space with a measure $\alpha$, and $h:A\to
B$ is a measurable map. The push-forward  of the measure $\alpha$
is the measure $\beta$ on $B$ defined by
$\beta(S)=\alpha(h^{-1}(S))$, where $S\subset B$.}
 of the
$K$-invariant measure $\mu$ under the map (\ref{radial-map}) is
given by the formula
\begin{equation}
 \mathrm{const}(\kappa)\cdot \prod_{1\le k<l\le n}
\bigl|e^{i\phi_k}-e^{i\phi_l}\bigr|^{ 2\kappa} \prod_{j=1}^n
d\phi_j \label{radial-part}
,\end{equation}
where

--- $\kappa=1/2$ for $\U(n)/\OO(n)$

--- $\kappa=1$ for $\U(n)\times\U(n)/\U(n)$

--- $\kappa=2$ for $\U(2n)/\Sp(n)$
\end{lemma}

See, for instance, \cite{Hel1}, X.1, various calculations of this
type are contained in \cite{Hua}.

\smallskip

Consider the space $L^2(K/H)^H$ consisting of $H$-invariant
$L^2$-functions. To each function $F\in L^2(K/H)^H$, we consider
the corresponding function $f$ in variables $\phi_j$. By Lemma
3.1, we have
$$
\int_{K/H} F_1(z)\ov{F_2(z)}d\mu(z)= \mathrm{const}
\int_{\T^n}f_1(\phi)\ov{f_2(\phi)}
 \prod_{1\le k<l\le n}
\bigl|e^{i\phi_k}-e^{i\phi_l}\bigr|^{ 2\kappa} \prod_{j=1}^n
d\phi_j \label{radial-product}
$$
This is precisely the inner product (\ref{jack-polynomials}).

\smallskip

 {\bf \punct Spherical functions. Preliminaries.}
Consider a symmetric
 space
$$K/H=\U(n)\times \U(n)/\U(n),\quad \U(n)/\OO(n), \quad
\U(2n)/\Sp(n)
$$
The space $L^2(K/H)$ is a multiplicity-free direct sum of
(finite-dimensional) $K$-modules $V_\nu$,
\begin{equation}L^2(K/H)\simeq\oplus_\nu V_\nu
\label{l2}
\end{equation}
 An irreducible $K$-module
$V_\nu$ participates in this direct sum iff  $V_\nu$ contains a
nonzero $H$-invariant vector $v_\nu$%
\footnote{ This statement is a variant of the Frobenius
reciprocity, see \cite{Kir}, 13.1, 13.5.
 For
description of modules satisfying this property, see, for
instance, \cite{Hel2}.}.
 By Gelfand's theorem (see \cite{Hel2}),
an $H$-invariant vector $v_\nu\in V_\nu$ is unique up to a scalar
factor.

 The {\it spherical function} of a module $V_\nu$
 is function on $K$ defined by
  \begin{equation}
  \xi_\nu(k)  :=\langle k\cdot v_\nu,v_\nu\rangle_{V_\nu}
  \end{equation}
where $v_\nu$ is normalized by the condition $\|v_\nu\|=1$.
Obviously, for $h_1$, $h_2\in H$ we have
 $$
  \xi_\nu(h_1kh_2)=  \xi_\nu(k)
$$
This allows to consider the function $\xi_\nu$ as a function on
our symmetric space $G/H$ or as a function on the double coset
space $H\setminus G/H$.

 In expansion (\ref{l2}), we have $\xi_\nu\in V_\nu$. Hence, the
functions $\xi_\nu$ form an orthogonal basis in $L^2(K/H)^H$.

The following fact is well-known, see, for instance, \cite{Mac},
Chapter 7.

{\it $H$-spherical functions on $K/H$ are  the Jack polynomials}
(modulo some normalizing scalar factors.)
 the parameter
$\kappa$ was determined in Lemma 3.1.

\smallskip

{\bf \punct Reduction of our problem.} Thus we intend to
investigate a positivity of the inner product
$$
\langle F_1,F_2\rangle_{\sigma,\tau} =\iint_{K/H\times K/H}
L_{\sigma,\tau}(z,u)F_1(z)\ov{F_2(u)}\,d\mu(z)\,d\mu(u)
$$
where $L_{\sigma,\tau}=\det(1-zu^*)^{\{\sigma\|\tau\}}$ is the
distribution defined above.

Let $\lambda$ be a collection (\ref{signatures}),
$P_\lambda(x;\kappa)$ be the corresponding Jack polynomial, and
$V_\lambda\subset C^\infty(K/H)$ be the $K$-invariant subspace
containing $P_\lambda(x,\kappa)$.

\begin{lemma}
For $F\in C^\infty(K/H)$ consider its expansion
 $F=\sum_\lambda F_\lambda$, where $F_\lambda\in V_\lambda$.
 Then
 $$
 \langle F,G\rangle_{\sigma,\tau}=
 \sum_\lambda c_\lambda(\sigma,\tau)
 \int_{K/H} F_\lambda(z) \ov{G_\lambda(z)} d\mu(z)
 $$
where $c_\lambda(\sigma,\tau)$ are some constants.
\end{lemma}

{\sc Proof.} Consider the integral operator
$$
U_{\sigma,\tau} F(z)=\int_{K/H} L_{\sigma,\tau}(z,u) F(u)\,d\mu(u)
$$
Since the $K$-invariance of the kernel $L_{\sigma,\tau}$, the
operator $U_{\sigma,\tau}$ is an intertwining operator for $K$.
 Since the action of $K$ in $C^\infty$
 is multiplicity-free, the restriction of $U_{\sigma,\tau}$ to
$V_\lambda$ is a scalar operator (see \cite{Kir}), 8.3), i.e.,
\begin{equation}
U_{\sigma,\tau} F_\lambda(z)= \int_{K/H} L_{\sigma,\tau}(z,u)
F_\lambda(u)\,d\mu(u)=c_\lambda(\sigma,\tau)F_\lambda (z)
\label{7}
\end{equation}
and this implies the required statement.
 \hfill $\square$

\smallskip

Next, we must evaluate the constants $c_\lambda(\sigma,\tau)$. For
this purpose, we substitute $F_\lambda=P_\lambda$ and

--- $z=1$ in the case $K/H=\U(n)$,

--- $z=1$ in the case $K/H=\U(n)/\OO(n)$,

--- $z=J$ (see (\ref{matrix-J})) in the case $K/H=\U(2n)/\Sp(n)$.

In all the cases we obtain the integral
\begin{multline}
\cL_\lambda(\kappa;\sigma,\tau)
=\\=
 \int_{\T^n} \prod_{k=1}^n
(1-e^{i\phi_k})^{\{\sigma\|\tau\}}
P_\lambda\bigl(e^{i\phi_1},\dots,e^{i\phi_n};\kappa\bigr)
 \prod_{1\le k<l\le n}
\bigl|e^{i\phi_k}-e^{i\phi_l}\bigr|^{ 2\kappa} \prod_{j=1}^n
d\phi_j \label{pseudo-kadell}
 \end{multline}
 where $\kappa$  is the same as in Lemma
3.1.

In the first two cases, we obtain this integral immediately, for
the case $\U(2n)/\Sp(n)$ we must write $\det(1-J\Omega)$, where
$J$ is (\ref{matrix-J}), and $\Omega$ is (\ref{can-form-usp}). The
matrix $\det(1-J\Omega)$ is diagonal, and
$$\det(1-J\Omega)=\prod_{k=1}^n (1-e^{i\phi_k/2})^{\{\sigma\|\tau\}}
(1-e^{-i\phi_k/2})^{\{\sigma\|\tau\}}=
\prod_{k=1}^n(1-e^{i\phi_k})^{\{\sigma\|\tau\}}
$$

{\bf \punct The Kadell integral.} Let $\kappa$ be a positive
integer. Let $\lambda$ satisfies (\ref{young}).
 The
Kadell integral (see \cite{Kad}, \cite{Mac}, VI.10, Example 7) is
given by
\begin{multline} \cK_\lambda(\kappa; r,s):=
\frac 1{n!} \int\limits_{[0,1]^n}
 P_\lambda(x;\kappa) \prod_{k=1}^n x_k^{r-1}(1-x_k)^{s-1}
 \prod_{1\le k<l\le n} |x_k-x_l|^{2\kappa}
\prod_{j=1}^n dx_j =\\=a_\lambda(\kappa;r,s) v_\lambda(\kappa)
\label{kadell}
\end{multline}
 where
\begin{align}
v_\lambda(\kappa)&= \prod_{1\le k<l\le n}
\frac{\Gamma(\lambda_k-\lambda_l+\kappa(l-k+1))}
{\Gamma(\lambda_k-\lambda_l+\kappa(l-k))} \label{v-lambda}
\\
a_\lambda(\kappa;r,s)&=\prod_{j=1}^n
\frac{\Gamma(\lambda_j+r+\kappa(n-j))\Gamma(s+\kappa(n-j))}
{\Gamma(\lambda_j+r+s+\kappa(2n-j+1))} \label{a-lambda}
\end{align}

{\bf \punct Transformation of the Kadell integral.}

\begin{prop}
Let $\Re(\sigma+\tau)>-1$, $\Re \kappa>0$, and $\lambda$ satisfies
(\ref{signatures}).
 Then the integral
$\cL_\lambda(\kappa;\sigma,\tau)$ given by (\ref{pseudo-kadell})
equals
\begin{multline}
\cL_\lambda(\kappa;\sigma,\tau)=
 (2\pi)^n n!\, v_\lambda(\kappa)
\times \\ \times
  \prod_{j=1}^n
\frac{(-1)^{\lambda_j}\Gamma(\sigma+\tau+1+\kappa(n-j))}
{\Gamma(-\lambda_j+\tau+1+\kappa(j-1))\Gamma(\lambda_j+\sigma+1+\kappa(n-j))}
\label{final-formula}
\end{multline}
 where $v_\lambda(\kappa)$ is the same as
above (\ref{v-lambda}).
\end{prop}

{\sc Proof.} First, let $\kappa$ be a positive integer. Assuming
$x_k=e^{i\phi_k}$,
$$\bigl|e^{i\phi_k}-e^{i\phi_l}\bigr|^2=
|x_k-x_l|^2=(x_k-x_l)(x_k^{-1}-x_l^{-1})=-(x_k-x_l)^2 x_k^{-1}
x_l^{-1}
$$
Next, transform the expression (\ref{pseudo-kadell}) to a contour
integral
\begin{multline*} \cL_\lambda(\kappa;\sigma,\tau)=
e^{in\pi \tau } (-1)^{n(n-1)\kappa/2} i^{-n}\times\\ \times
\int\limits_{|x_1|=\dots=|x_n|=1}
 P_\lambda(x;\kappa) \prod_{k=1}^n x_k^{-\tau -\kappa(n-1)-1}(1-x_k)^{\sigma+\tau}
 \prod_{1\le k<l\le n} (x_k-x_l)^{2\kappa}
\prod_{j=1}^n dx_j
\end{multline*}
Now the integrand is holomorphic in the polydisc $|z_k|<1$ outside
the cut $x_1\in [0,1]$, \dots, $x_n\in [0,1]$. We define the
branches of the factors by
$$
(1-x)^\mu\Bigr|_{x=0}=1,\qquad x^\nu\Bigr|_{x=e^{i\phi},\, \phi\in
[0,2\pi)}=e^{i\nu\phi}
$$
The domain of convergence of this integral is
\begin{equation}
\Re(\sigma+\tau)>-1 \label{domain1}
\end{equation}

 Then we deform each one-dimensional contour $x_j=e^{i\phi_j}$
to the contour around the segment $x_j=[0,1]$ lying in an
infinitely thin strip $|\Im x_j|<\epsilon$. Thus our expression
converts to
$$
 (-1)^{n(n-1)\kappa/2}(-2\sin \pi \tau)^n \int_{[0,1]^n}
 \Bigl\{\text{the same integrand}\Bigr\}
 $$
The domain of convergence is reduced to
\begin{equation}
\Re(\sigma+\tau)>-1, \qquad \Re\tau<-\kappa(n-1) \label{domain2}
 \end{equation}

 Applying the Kadell formula, we obtain
  \begin{multline}
  \cL_\lambda(\kappa;\sigma,\tau)
=2^n n! (-\sin(\tau\pi))^n (-1)^{n(n-1)\kappa/2}
 v_\lambda(\kappa)\times\\ \times
\prod_{j=1}^n
\frac{\Gamma(\lambda_j-\kappa(j-1))\Gamma(\sigma+\tau+1+\kappa(n-j))}
{\Gamma(\lambda_j+\sigma+1+\kappa(n-j))}
\label{expression-for-positivity}
\end{multline}
 Then we write
\begin{multline*}
-\sin (\pi\tau)\Gamma(\lambda_j-\tau-\kappa(j-1))=
\\=
(-1)^{\lambda_j-\tau-\kappa(j-1)}
\sin(\lambda_j-\tau-\kappa(j-1))\pi)
\Gamma(\lambda_j-\tau-\kappa(j-1))
\end{multline*}
(again, we use the condition $\kappa\in \Z$). Applying the
reflection formula
$$
{\sin(\pi z)}\Gamma(z)= \pi/{\Gamma(1-z)}
$$
we obtain (\ref{final-formula}).

Thus, we have the required identity (\ref{final-formula}) under
the condition (\ref{domain2}). But the both sides of
(\ref{final-formula}) are analytic in the domain (\ref{domain1}).
Hence the identity is valid in this domain.

\smallskip

Thus the statement is proved for an integer positive $\kappa$.

\smallskip

Next, we intend to apply the following {\it Carlson theorem}, see
\cite{AAR}, 2.8.1.

\smallskip

--- {\it Let $f(z)$ be an analytical function in the domain $\Re
z>0$ and $f(z)=O(e^{\alpha|z|})$ with some $\alpha<\pi$. If
$f(n)=0$ for $n=1$, $2$,\dots, then $f(z)$ is identically zero.}

\smallskip

Fix real positive $\sigma$, $\tau$. The identity
(\ref{pseudo-kadell}) is valid for positive integer $\kappa$. We
intend to show that it is valid for $\Re\kappa>0$. First, we
slightly change the factors of the integrand
 $$
\bigl|e^{i\phi_k}-e^{i\phi_l}\bigr| \mapsto
\bigl|\tfrac12(e^{i\phi_k}-e^{i\phi_l})\bigr| ,\qquad
 (1-e^{i\phi_k})
\mapsto \tfrac12(1-e^{i\phi_k})
$$
 Accordingly we multiply the right-hand side of
(\ref{pseudo-kadell}) by $2^{-n(n-1)\kappa-n(\sigma+\tau)}$. Now
the expression $\prod_{j=1}^n\dots\prod_{1\le j<k\le n}\dots$ in
the integrand is $<1$.

 For coefficients $u_\mu$ of the Jack polynomials in formula
(\ref{u}), there exists a semi-explicit  expression that is
rational in the variable $\kappa$, see \cite{Mac}, VI.10, page
379. The poles of these expressions are some non-positive rational
numbers. Hence, for a fixed $\lambda$, the expression
$P_\lambda(x;\kappa)$ admits a holomorphic continuation to the
domain $\Re\kappa>0$, and moreover $P_\lambda(x;\kappa)$ has a at
most a polynomial growth in $\kappa$.

Thus our integrand has (at most) a polynomial growth in $\kappa$.
Hence the same statement is valid for the integral.

\smallskip

Now we estimate the growth of the right-hand side of
\ref{final-formula}. First,  our formula is valid for
$\lambda_1=\dots=\lambda_n=0$,  (in this case, we obtain   the
Cauchy-type form of the Selberg integral, see \cite{AAR}; see also
\cite{Mac}, (10.38),  for $\lambda=0$). The integrand now is less
than 1, and-hence the right-hand side also is bounded.

 Consider the ratio of right-hand
sides
$$\cL_\lambda(\kappa;\sigma,\tau)/\cL_0(\kappa;\sigma,\tau)$$
This ratio is the product of expressions having the form
$\Gamma(\mu+l\kappa)/\Gamma(\nu+l\kappa)$. This product has a
polynomial growth in $\kappa$, since
$$
\frac{\Gamma(a+z)}{\Gamma(b+z)} \sim z^{a-b},\qquad |z|\to\infty
$$
(see, for instance \cite{AAR}, 1.4.3; the
 asymptotic is valid in
the domain $|\arg z|<\pi-\epsilon$).

Now we have two functions having a polynomial growth as
$|z|\to\infty$, they coincide in positive integers and  by the
Carlson theorem they coincide in the half-plane.

\smallskip

It remains to omit the condition $\lambda_n\ge 0$. But the
transformation
$$\lambda_j\mapsto\lambda_j-1,\qquad \tau\mapsto \tau-1,\qquad
\sigma\mapsto\sigma+1
$$
does not change the both sides of our integral, and hence our
formula can be used for general $\lambda$ satisfying
(\ref{signatures})

{\bf \punct  Positivity.} We must analyze, when
(\ref{final-formula}) has constant signs for all $\lambda$
satisfying (\ref{signatures}). The factor $v_\lambda$ given by
(\ref{v-lambda}) is positive. Hence we must analyze positivity of
the factor $\prod\dots$ in (\ref{final-formula}). Equivalently we
must trace the constancy of the sign of the factor
$$\prod_{j=1}^n
\frac{\Gamma(\lambda_j-\tau-\kappa(j-1))}
{\Gamma(\lambda_j+\sigma+1+\kappa(n-j))}
$$
  in
(\ref{expression-for-positivity}).

{\it The case $\kappa=2$.} We must analyze positivity of the
expression
\begin{equation}
\prod_{j=1}^n
 \frac{\Gamma(\lambda_j-\tau-2j+2)}
{\Gamma(\lambda_j+\sigma+1 +2(n-j))} \label{product-for-signs}
\end{equation}
 The function $\Gamma(x)$ changes its
sign at the points $x=0$, $-1$, $-2$, \dots. Hence the condition
\begin{equation}
[-\tau-2j+2]= [\sigma+1+2(n-j)]\qquad \Longleftrightarrow\qquad
[-\tau]=[\sigma+2n-1] \label{8}
\end{equation}
 is sufficient for
positivity.

Let us show that this condition is necessary. Assume that all the
$\lambda_j$ are negative integers having large absolute values,
and assume that $\lambda_i-\lambda_j$ also are large.
 Then the
sign of (\ref{product-for-signs}) depends only on parities of
$\lambda_j$. Also, it is clear, that all the factors in
(\ref{product-for-signs}) have to be positive. But this implies
(\ref{8}).

{\it The case $\kappa=1$} is similar. We must examine the
positivity of
$$
\prod_{j=1}^n \frac{\Gamma(\lambda_j-\tau-j+1)}
{\Gamma(\lambda_j+\sigma+n-j+1)}
$$

{\it The case $\kappa=1/2$.} We must examine the positivity of
$$
\prod_{j=1}^n \frac{\Gamma(\lambda_j-\tau-j/2+1/2)}
{\Gamma(\lambda_j+\sigma+1+(n-j)/2)}
$$
Considering even $j$-s, we obtain
$$[-\tau+1/2]=[\sigma+1+n/2]$$
Considering odd $j$-s, we obtain
$$[-\tau+1/2+1/2]=[\sigma+1+n/2+1/2]$$
This implies Theorem 2.5.

\bigskip

{\bf\large 4. Remarks on general case}

\medskip

\addtocounter{sec}{1} \setcounter{equation}{0}
\setcounter{punct}{0} \setcounter{theorem}{0}

Here we discuss Stein--Sahi representations
of arbitrary classical groups.

\smallskip

{\bf \punct Flat matrix spaces.} We consider the following 10
types of  spaces $\cM$ of $n\times n$ matrices;
in each case we write 
  groups $G\supset K\supset H$, the sence of this
  notation
  will be explained below; if $n$ have to be odd, then
  we write $n=2k$

\medskip

---$n\times n$ matrices over $\C$;

\hfill $\GL(2n,\C)\supset\U(2n)\supset \U(n)\times\U(n)$;

\smallskip

--- $n\times n$ Hermitian matrices  over $\C$;
 
\hfill $\U(n,n)\supset \U(n)\times\U(n)\supset \U(n)$

\smallskip

--- $n\times n$ symmetric matrices  over $\C$;

\hfill $\Sp(2n,\C)\supset\Sp(n)\supset\U(n)$

\medskip

--- $2k\times 2k$ skew-symmetric matrices over $\C$;

\hfill  $\SO(2k,\C)\supset\SO(2k)\supset\U(k)$.

\smallskip

--- $n\times n$   matrices over $\R$;
 
 \hfill $\GL(2n,\R)\supset
\OO(2n)\supset \OO(n)\times\OO(n)$

\smallskip

--- $n\times n$ symmetric matrices over $\R$;
 
 \hfill $\Sp(2n,\R)\supset
\U(n)\supset\OO(n)$

\smallskip

--- $2k\times 2k$ skew-symmetric matrices over $\R$

\hfill
$\SO(2k,2k)\supset \OO(2k)\times\OO(2k)\supset\O(2k)$

\medskip

--- $n\times n$ matrices over $\H$,

\hfill$\GL(2n,\H)\supset \Sp(2n)\supset\Sp(n)\times \Sp(n)$

\smallskip

--- $n\times n$ anti-Hermitian matrices over $\H$,

\hfill$\Sp(n,n)\supset \Sp(n)\times\Sp(n)\supset \Sp(n)$

\smallskip

--- $k\times k$ Hermitian matrices over $\H$;
 
 \hfill $\SOS(4k)\supset \U(2k)\supset \Sp(k)$

\medskip

The group $G$ is isomorphic to the group of linear-fractional
transformations
$$x\mapsto x^{[g]}= (a+xc)^{-1}(b+xd)$$
 preserving the symmetry
condition. The Jacobian of such transformation is
$$
J(g,x)=\det|(a+xc)|^{-h}
\Bigl|\det\begin{pmatrix}a&b\\c&d\end{pmatrix}\Bigr|^{\beta}
$$
where
$$h=\frac 2n \dim\cM,$$
$n$ is  size of the matrices,
and $\beta=n\dim\K$ in 3 cases then we consider all the matrices
without conditions of symmetry (i.e., $G=\GL(n,\K)$); in all the
remaining cases
$\bigl|\det\begin{pmatrix}a&b\\c&d\end{pmatrix}\bigr|=1$, see
\cite{Ner-beta}, Lemma 1.3.

 The
group $K$ is the maximal compact subgroup  in $G$, and $H\subset
K$ is the stabilizer of $0$. In all the cases, $K/H$ is a compact
symmetric space, and $\cM$ is an open dense subset in $K/H$

{\sc Remark.} In all the cases, 
there exists  some Grassmannian type
realisation of $K/H$ in spirit of our Subsections 1.3, 2.2, 2.5,
see a  table 
in \cite{Ner-krein}).

We define the Stein--Sahi inner product in $C^\infty(\cM)$ by
\begin{equation}
\langle F_1,F_2\rangle =\iint_{\cM\times\cM} |\det(x-y)|^\theta
F_1(x)\ov {F_2(y)} dx\,dy
\label{general-product}
\end{equation}
The group $G$ acts  in the space of functions on $\cM$ by formula
\begin{equation}
\rho\begin{pmatrix} a&b\\c&d\end{pmatrix}F(x)
=F\bigl((a+xc)^{-1}(b+xd)\bigr)|\det(a+xc)^{-h-\theta}
\label{general-action}
\end{equation}
these transformations preserve the Hermitian form $\theta$.

\smallskip

{\bf \punct Existence of an island of positivity.}
 First, consider
the representation of $G$ in functions on $\cM$ given by
\begin{equation}
F(x)\mapsto f\bigl((a+xc)^{-1}(b+xd)\bigr) \det(a+xc)^{-h/2}
\label{principal-series}
 \end{equation}
 It is a representation of a
degenerated unitary principal series.

Consider the maximal parabolic subgroup $P\subset G$ consisting of
matrices having the form $\begin{pmatrix} a&b\\0&d\end{pmatrix}$.
This subgroup acts on functions by affine transformations
\begin{equation}
F(x)\mapsto F(a^{-1}(b+x d))\det (a)^u \det(d)^v \label{15}
,\end{equation}
where $u$, $d$ are some explicit real numbers
(non-interesting for us).

\begin{lemma}
If $G\ne \Sp(2n,\R)$, $\U(n,n)$, $\SOS(4n)$,  then the
representation (\ref{15}) of $P$ is irreducible.
\end{lemma}

{\sc Proof.} Consider the subgroup $N\subset P$ consisting of
matrices $\begin{pmatrix} 1&b\\0&1\end{pmatrix}$, it acts by
shifts $f(x)\mapsto f(x+b)$.

Also consider the subgroup $L\subset P$ consisting of  matrices
$\begin{pmatrix} a&0\\0&d\end{pmatrix}$, it acts by the
transformations $f(x)\mapsto f(axd^{-1})$.

After the Fourier transform
$$F(x)\mapsto \int_\cM F(x)\exp(i\Re\mathrm{tr}\, x\xi^*)\,dx$$
the group $N$ acts via multiplications by functions
\begin{equation}
\Phi(\xi)\mapsto \Phi(\xi)\exp (i\Re\mathrm{tr}\, \xi b^*) \label{16}
\end{equation}

Assume that the representation (\ref{15}) is reducible. Then there
exists an intertwining projector $\Pi$. In the Fourier model, it
commutes with all the operators (\ref{16}) and hence $\Pi$ is a
multiplication by a function $\chi$ taking values $0$ and $1$. But
this operator also must commute with the subgroup $L$, which acts by
the transformations having the form
$$\Phi(\xi)\mapsto \Phi(A \xi D^{-1})|\det A|^{u'}|\det D|^{v'}$$
But in all the cases, the group of transformations $\xi\mapsto
A\xi D^{-1}$ has an open orbit on $\cM$ and hence $\chi$
is a constant.
\hfill $\square$

\begin{cor}
If  $G\ne\U(n,n)$, $\Sp(2n,\R)$, $\SOS(4n)$, then the
representation (\ref{principal-series}) is irreducible.
\end{cor}

\smallskip

Since (\ref{principal-series}) is irreducible, it admits a unique
invariant Hermitian form and this form is the $L^2$-inner product.
Representations (\ref{general-action})
 lying near (\ref{principal-series}) admit
the invariant Hermitian form (\ref{general-product}), 
and since this Hermitian form is a minor
deformation of $L^2$-inner product, this Hermitian form also is
positive definite. It is not difficult to justify these heuristic
arguments in this case (there are abstract theorems of this type,
see \cite{KS}, Theorems 6, 8; unfortunately they do not cover our
cases).

Description of the  poles of the distribution $|\det(x)|^\theta$
that appear in (\ref{general-product})
 is a relatively standard problem, see
Sato, Shintani \cite{SAS}%
 \footnote{For more information on determinant distributions
 (explicit expressions at poles, supports, etc.)
 see works  Rais \cite{Rai}, 
  Ricci, Stein \cite{RS}, Muro \cite{Mur}.}. 
  Knowledge of these poles allows to find
precisely the interval of positivity. We omit details.

\smallskip

{\bf \punct Other possible ways of proof.} Existence of 
Stein--Sahi
series is a relatively simple fact, and apparently it can be proved uniformly
in various  ways.

One possibility is to use the work of Branson, Olafsson,
\O rsted, \cite{BOO}.

Another way is to use our method of Sections 1-3, but a necessary
identity for multivariate Jacobi polynomials (see \cite{HS}) 
in my knowledge
 is not  known yet. But for the groups
$G=\OO(2n,2n)$, $\Sp(n,n)$, $\GL(2n,\C)$ this 
way can easily be
realized (the first two cases 
are discussed in \cite{Ner-sobolev})
due existence of explicit elementary formulas
(Weyl's and Berezin--Karpelivich's) for spherical
functions.

\bigskip

{\bf \large 5. One problem of non-$L^2$ harmonic analysis.}

\addtocounter{sec}{1} \setcounter{equation}{0}
\setcounter{punct}{0} \setcounter{theorem}{0}

\medskip

{\bf  A. Construction of kernels}

{\bf \punct Uniform
realizations of symmetric spaces.}
 Consider the following
classical groups
\begin{multline}
Q=\GL(n,\C),\,\GL(n,\R),\,\GL(n,\H),\,\OO(p,q),\,\U(p,q),\,\Sp(p,q),\,
\OO(n,\C),\,
\\
\Sp(2n,\C),\, \Sp(2n,\R),\,\SOS(2n) \label{list-of-groups}
\end{multline}
  
  {\sc Remark.} Our list contains $\GL(n,\K)$ and
not $\SL(n,\K)$, also $\U(p,q)$ and not $\SU(p,q)$. Modulo this
remark, our list contains all the classical groups.

\smallskip

 Consider an affine symmetric
space having the form $Q/Y$, where $Y$ is a symmetric subgroup in
$Q$  (according the Berger classification, there are 54 series of
classical affine symmetric spaces).
 It turns out to be (see
\cite{Ner-uniform}) that the space $Q/S$ admits a
realization as a set, whose points are pairs of complementary
subspaces $(V,W)$ in some linear space $\K^m$, these subspaces
 satisfy  some simple conditions (as an isotropy with
respect to some form, orthogonality, or existing of a given
involution transposing two subspaces).

\smallskip

{\sc Example 1.} {\it The unitary group $\U(n)$.} 
Consider the
space $\C^n\oplus\C^n$ equipped with the Hermitian form $H$
having the matrix 
$\begin{pmatrix} 1&0\\0&-1\end{pmatrix}$ as above (Subsection 1.2). 
Fix the linear operator $J$ in
$\C^n\oplus\C^n$ having the matrix $\begin{pmatrix}
1&0\\0&-1\end{pmatrix}$. Consider the set $\mathcal U$,
 whose
points are pairs of subspaces $(V,W)$ such that

$0^*$. $\C^n\oplus\C^n=V\oplus W$

$1^*$. $V$ and $W$ are $H$-isotropic.

$2^*$. $W=JV$

Now let $z:\C^n\to\C^n$ be an unitary operator. Let $V$ be a graph
of $z$, and $W$ be the graph of $(-z)$. Then $V$, $W$ satisfy to
the conditions 1-3. It can be readily checked that all the pairs
satisfying 1-3 have this form.

Thus $\mathcal U\simeq\U(n)$.%
\footnote{ Of course, the subspace $W$ in our construction seems artificial
 and can be forgotten (as it was done above in Sections 1-2).
 But it take part in a general construction
 of the kernel below).}

\smallskip

{\sc Example 2.} {\it The space $\U(n,n)/\GL(n,\C)$.} Consider the
same space $\C^n\oplus\C^n$ equipped with the same Hermitian form
$H$. Consider the set $\mathcal H$, whose points are pairs $V$,
$W$ of subspaces satisfying the conditions $0^*-1^*$ from the
previous example. Obviously, $\mathcal H\simeq \U(n,n)/\GL(n,\C)$.

\smallskip

{\sc Example 3.} {\it The spaces
$\U(n,n)/\U(p,q)\times\U(n-p,n-q)$.} Consider the same space
$\C^n\oplus\C^n$ with the same form $H$. Consider the set
$\mathcal G$, whose points are pairs of subspaces $(V,W)$
satisfying the conditions $0^*$  and

--- $\dim V=\alpha$

--- $W$ is the orthocomplement of $V$ with respect to $H$.

Open orbits of the group $\U(n,n)$ 
on $\mathcal G$ are enumerated by
the inertia indices of $H$ on $V$. Thus%
\footnote{Again, it seems that $W$ is an artificial element of the
construction (since $W$ is determined by $V$), but the
transversality condition $0^*$ involves $W$.}
$$\mathcal
G\simeq\cup_{p+q=\alpha}
\U(n,n)/\U(p,q)\times\U(n-p,n-q).$$

It appears that enumerating carefully all the possible conditions
of this type, we obtain precisely the classical part of the Berger
list.%
\footnote{We must consider all the fields $\K=\R$, $\C$, $\H$,
all  possible natural forms and all possible involutions.
For a table, see \cite{Ner-uniform}.}

\smallskip

{\bf \punct Overgroup of the symmetric space.}
In Examples 1 and 3, the subspace $W$ is determined by the
subspace $V$, and hence we have an embedding of the symmetric
space to some Grassmannian  $G/P$%
\footnote{In both examples $G=\U(n,n)$}; 
the same phenomenon
holds for 44 series of symmetric spaces. 

In Example 2, the
subspaces $V$, $W$ are 'independent', and we obtain an open
embedding of a symmetric space to a product of two Grassmannians
$G/P_1\times
G/P_2$%
\footnote{In our example $G=\U(n,n)$}
 (and this happens for remaining
10 series).

In particular, we obtain that {\it for each classical symmetric
space
$Q/Y$ there exists%
\footnote{Emphasis that $Q$ is contained in the list
(\ref{list-of-groups}). For instance, $\U(n,n)$ acts on the space
$\U(n)$; but there is no  group $G\supset\SU(n)\times\SU(n)$
acting on the hypersurface $\SU(n)\subset\U(n)$}
 larger group $G$
acting locally on $Q/Y$}.

\smallskip

{\sc Remark.} It seems that this (trivial) observation was firstly
claimed in \cite{Ner-uniform}. But overgroups of symmetric spaces
were discussed by different reasons by many authors, in particular
Hua Loo Keng, Nagano, Takeuchi, Goncharov, Gindikin, Kaneyuki.
Nagano \cite{Nag} gave a complete analysis of such
overgroups in Riemannian case (including exceptional cases).
Makarevich  \cite{Mak} in 1973 classified open  orbits of
reductive groups in Grasmannans. In fact, the right-hand side of
long Makarevich's  tables includes the whole classical part of
Berger's list, and this implies our observation.
 As far as I know,
nobody actually  compared these Bergers's 
and Makarevich's lists during 25 years.%
\footnote{Only few exceptional symmetric spaces admit such
overgroup;  see
\cite{Bert}.}

\smallskip

{\bf \punct Double ratio of 4 subspaces.} Consider a linear space
$\K^{p+q}$. Let $V_1$, $W_1$, $V_2$, $W_2$ be four subspaces in
$\K^{p+q}$, and
$$
\dim V_1=\dim V_2=p,\qquad\dim W_1=\dim W_2=q
$$
For a quadruple being in {\it general position}, we define a
canonical operator ({\it the Hua Loo Keng double ratio})
$\frD:V_1\to V_1$
 by the following rule. Since $\K^{p+q}=V_1\oplus W_1$, the
subspace   $V_2$ is a graph of an operator $A:V_1\to W_1$, and $B$
is a graph of an operator $B:W_1\to V_1$. We assume
$$\frD(V_1,W_1,V_2,W_2)=BA$$

The operator is canonically defined, in particular, its
eigenvalues are invariants of a quadruple of subspaces.

By Subsection 5.1, the double ratio is well defined in
each classical symmetric space.

\smallskip

 {\bf \punct A problem.
  Different formulations.} The author think, that
  for the
following problems of the harmonic analysis,
the explicit Plancherel formula can be obtained.
 The arguments for
support of this point of view are discussed below.

We use the notation of Section 4.

\smallskip

a) {\it Tensor products.}
 For a given group $G$ from the list
 (\ref{line1}--{line2}), decompose a tensor
 product of two Stein--Sahi representations.

\smallskip

b) {\it Restrictions.} Let $Q/Y$ be a symmetric space.
Assume $G$ that its overgroup has the Stein--Sahi
representations. Assume that the local action
of $G$ on $Q/Y$ is locally isomorphic 
to the action $G$ on the matrix space, see 4.1.
Decompose the restriction
of a Stein--Sahi representation
to $G$.

\smallskip

{\sc Remark.} The problem about tensor products is
 a  problem of the restriction of a representation $\rho\otimes
\rho'$ of the group $G\times G$ to the diagonal subgroup $G$.

\smallskip

c) {\it Natural kernels on symmetric spaces.} Let $Q/Y$ be
 the same
as above. Let restrict the Stein kernel $L(\cdot,\cdot)$ to $Q/Y$.
Consider the inner product in $C^\infty(Q/Y)$ defined by the
kernel $L$ (see (\ref{general-product})) and the unitary representation of $Q$ in this
space. To find the Plancherel formula for this representation.

\smallskip

{\sc Remark.} This question slightly differs from the previous one
in the case, then $Q$ has several open orbits $Q/Y_j$ on
the Grassmannian,
see Example 3 in 5.1.

\smallskip

d) {\it More general kernels on symmetric spaces.} Consider a
classical symmetric space $Q/Y$ realized as in 5.1, i.e., points
of $Q/Y$ are pairs $r=[V,W]$ of complementary subspaces. We define
(see Addendum to \cite{Ner-determinant}) the kernel
$K(\cdot,\cdot)$ on $Q/Y$ by
$$
K([V_1,W_1], [V_2,W_2])= \Bigl|\det\frac{\frD(V_1,W_1,V_2,W_2)}
{1-\frD(V_1,W_1,V_2,W_2)}\Bigl|^\theta
$$
and the inner product on $C^\infty(Q/Y)$ given by
$$
\langle f_1,f_2\rangle= \iint_{Q/Y\times Q/Y}
K(r_1,r_2)f_1(r_1)\ov{f_2(r_2)} \,d\mu(r_1)\, d\mu(r_2)
$$
where $\mu$ is the $Q$-invariant measure on $Q/Y$.

If the overgroup $G$ of $Q/Y$ admits Stein--Sahi representations%
\footnote{To avoid ambiguity, we also must require that the local
action of $G$ on $Q/Y$ is equivalent to the action on the matrix
spaces (or Grassmannians) from 4.1. This slip in  speech is
important only for $G=\GL(n,\K)$.}, then our construction is
equivalent to the previous one.

Otherwise, we obtain a problem of positivity of the inner product.
There are spaces $Q/Y$, for which islands of positivity are absent
(\cite{Ner-nonexisting}), 
and the cases, when such islands exist
(\cite{Mol-sl3}, \cite{Ner-nonexisting}). Generally, the problem
is open.

Again, our question  is to find the Plancherel measure.

I do not know has this question sense if there is no positivity,
see a discussion in Addendum to \cite{Ner-determinant}.

\smallskip

e) {\it More tensor products.} Our Stein representations of $G$
are induced from some parabolic $P$. 
We ask about the tensor
products of arbitrary two representations induced from this
parabolic.

\smallskip

{\bf B. Discussion}

\nopagebreak

\smallskip

Below we discuss the variant b) of the problem, i.e., the
restriction of a Stein representation to a symmetric subgroup.

\smallskip

{\bf \punct Approximation of $L^2$ on symmetric spaces.} First,
for all the groups $G$,  Stein--Sahi series includes  some
representation (or representations) of a principal series (in the
notations of 1.8, 2.3, 2.7, this corresponds to $s=t$, in the
notation of section 4, this corresponds to $\theta=-h/2$.

In this case,
 our problem is equivalent to decomposition of $L^2(Q/Y)$.

Hence our representations are some kind of deformation of
$L^2(Q/Y)$

\smallskip

{\bf \punct Modern picture of $L^2$ on a pseudo-Riemannian
symmetric space.} The problem of decomposition of $L^2$ on an
arbitrary pseudo-Riemannian symmetric space was considered as
important after the Flensted-Jensen's work  on
discrete series \cite{FJ}, 1980.

After large efforts, the problem is not
completely solved up to now.

First, in the following partial cases, the explicit solution
is known.

--- Riemannian symmetric spaces $G/K$ (Gelfand--Naimark,
\cite{GN2}
 for complex classical groups, Gindikin--Karpelevich,
\cite{GK2}, for general case).

--- Semisimple Lie groups (Gelfand--Naimark, \cite{GN2} for complex
classical groups, Harish-Chandara, \cite{HC}, for general case

--- Rank 1 spaces, Molchanov (\cite{Mol-rank1},
\cite{Mol-VINITI}).

--- The spaces of the form $G_\C/G_\R$, where $G_\C$ is a complex
group, and $G_\R$ is a real form (Harinck, \cite{Har1},
\cite{Har2}).

\smallskip

Second, there was a long story of search of the general Plancherel
formula. Author do not intend discuss it and fix modern situation.

There are
 two groups of authors, van den Ban -- Schlichtkrull
\cite{Ban1}, \cite{Ban2} and Delorme -- Carmona \cite{Del}
proposed variants of a Plancherel formula
(earlier a Plancherel formula was announced by
Oshima, but proof was not published). Proofs are very heavy
and  occupy  large collections of papers; also
these formulae
are nonexplicit and the $c$-function is not evaluated.

Oshima \cite{Osh} recently announced (with sketches of proofs) an
explicit formula for $c$-function for the "most continuous part of
spectrum". Apparently, Oshima's formula together with van den
Ban--Schlichtkrull's residue calculus allow to receive a final
solution.

A more important problem is to make this
branch of analysis available to a wider
mathematical community. 

\smallskip

{\bf\punct Degenerated cases of our problem.}
 Consider the
representations $\rho_{\sigma,\tau}$ of $\U(n,n)$ defined  in 1.6.
Consider
 a tensor product $\rho_{\sigma,\tau}\otimes\rho_{\sigma',\tau'}$.
 By definition, this representation acts in a space of functions
 on $\U(n)\times\U(n)$. The  diagonal group
 $\U(n,n)\subset\U(n,n)\times\U(n,n)$ has an open orbit
$\U(n,n)/\GL(n,\C)$ on $\U(n,n)$%
\footnote{Indeed, a point of $\U(n)$ can be considered as
a pair of maximal
isotropic subspaces in $\C^n\oplus\C^n$, see Example 2 from 5.1.
 A stabilizer of
a pair of maxomal isotropic subspaces is $\GL(n,\C)$.}.

Thus, $\U(n,n)$ acts in the space of functions
on $\U(n,n)/\GL(n,\C)$.

 Representations
$\rho_{\sigma,0}$ are highest weight representations of $\U(n,n)$,
and $\rho_{0,\tau}$ are lowest weight representations. The
problems of decomposition
$$
\rho_{\sigma,0}\otimes \rho_{\sigma',0},\qquad
\rho_{\sigma,0}\otimes \rho_{0,\tau}
$$
are  degenerated cases of our problem.

In the first case, we have purely discrete spectrum, and problem
has combinatorial nature (see, for instance \cite{JV}), in the
second case we obtain 'Berezin representations' discussed in
Introduction.

Hence our problem admits a degeneration into
2 important problems.

\smallskip

Restriction of a highest weight representation $\rho_{\sigma,0}$
of $G$ to a symmetric subgroup $Q$  leads to a similar
alternative,
i.e., we have purely discrete spectrum%
\footnote{An example is  $G=\U(n,n)$, $Q=
\U(p,q)\times\U(n-p,n-q)$, the corresponding symmetric spaces are
described in Example 3 of 5.1.} or  Berezin representation%
\footnote{An example is $\U(n,n)\supset \OO(n,n)$. This
restriction problem can be interpreted as a problem of analysis on
the symmetric space $\OO(n,n)/\OO(n,\C)$.}

 But generally such restriction problems leads to harmonic
 analysis on some class of spaces that is more general than
 symmetric spaces.%
 \footnote{Some counterexamples are the  restriction problems
 $\U(p,q)\supset \OO(p,q)$ for $p\ne q$ and
 $\U(n,n)\supset\Sp(2n,\R)$ for odd $n$.}

\smallskip

{\bf\punct Rank 1 cases.}

\smallskip

a) {\it Tensor products of unitary representations
$\SL(2,R)=\SU(1,1)$.}  For {\it linear}
 representations of $\SL(2,\R)$, the spectra were determined by
 Pukanzsky \cite{Puk} and the Plancherel formula by
 Molchanov \cite{Mol-comp}, \cite{Mol-tensor}.

Even this relatively simple problem is unexpectedly
non-trivial
and not perfectly understood up to now%
\footnote{A unitary (projective) representation of $SU(1,1)$
depend on {\it two} real parameters, and a finite-dimensional
representation of $\SU(1,1)$ depends on {\it one} integer
parameter. By this reason, formulae for infinite-dimensional
representations are not analytic continuations of
formulae for finite-dimensional representations.}
 \footnote{This
problem has a continuous spectrum of multiplicity 2 and also some
discrete spectra. By this reason,
 it is a good proving ground for understanding of some
 general phenomena of harmonic analysis.}
 \footnote{Let $G$ be a semisimple group and $\rho_1$, $\rho_2$
 are its unitary representations. Usually the spectrum
 of the tensor product $\rho_1\otimes\rho_2$ has infinite
 multiplicities (non-formal argument: 'functional dimension' of
 this product is too large). In a general situation, the
 problem
 of description of specter can be reasonable, but
 question about the Plancherel formula looks
 as  dangerous.
 It seems that our problem of analysis
 on $\U(n,n)/\GL(n,\C)$ is the closest multirank
 approximation to tensor products for $\SU(1,1)$.}
, see recent works \cite{Ner-shift}, \cite{Ner-jacobi},
\cite{Gro}.

\smallskip

{\it Other rank one cases.} There are many explicit calculations
in the rank one case, that can be attributed to our subject, see
\cite{Nai1}-\cite{Nai2}, \cite{Muk}, \cite{Mol-sl3}, \cite{vDM},
our list do not pretend to be complete.

\smallskip

{\bf \punct Compact symmetric spaces.} Let our space $Q/Y$ is
compact. For instance, let $Q/Y$ is $\U(n)$. Then we have the
space of functions on $\U(n)$ with the inner product (1.15).
 The group
$\U(n)\times \U(n)$ acts in this space by transformations
$F(z)\mapsto F(a^{-1}zb)$. Search of the Plancherel formula is
equivalent to expansion of $L_{\sigma,\tau}$ in spherical
functions (modulo some correction factors).

Since we have obtained explicit expansions, 
we also know the Plancherel
formula for the corresponding compact symmetric spaces.

\smallskip

{\bf\punct Oshima's argument.} In \cite{Osh}, Oshima shows that
for two spaces $Q_1/Y_1$ and $Q_2/Y_2$ having the same
complexification,
 their $c$-functions satisfy the same difference equations.
In particular, their ratio a priory is a trigonometric function.

Apparently, this phenomenon survives in our case. Since in the
compact case the Plancherel formula can be obtained, it is natural
to believe that it can be obtained in a general case.%
\footnote{The author had evaluated the Plancherel measure in some
relatively simple cases, see \cite{Ner-nonexisting}.}

{\sf Math.Physics group,
 Institute of Theoretical and Experimental Physics, %
\linebreak
B.Cheremushkinskaya, 25, Moscow 117 259, Russia}


\tt neretin@mccme.ru

and

{\sf Math.Departnent,
University of Vienna,
Nordbergstrasse, 15,
Vienna, Austria}

\end{document}